\documentclass[12pt]{article}
\usepackage{fancyhdr,amsmath, psfrag, graphicx, amsfonts, layout, color,enumerate, eurosym, tikz,lscape}
\definecolor{forestgreen}{rgb}{0.13,0.54,0.13}
\usetikzlibrary{shapes}

\font\fontauthors=cmcsc10 scaled \magstep1

\def\Arg{\mathop{\rm Arg}\nolimits}
\def\be{\begin{equation}}

\def\Diag{\mathop{\rm Diag}\nolimits}

\def\ee{\end{equation}}
\def\egalenloi{\smash{\mathop{~= ~}\limits ^{\cal L}}}
\def\egalLoi{{~\mathop{= }\limits^{\rond L}}~}

\def\Proba{\g P}

\def\Sp{\mathop{\rm Sp}\nolimits}

\def\stirling2 #1#2{\left\{\begin{matrix} #1\\#2\end{matrix}\right\}}

\def\transp#1{{\vphantom{#1}}^t#1}

\def\Var{\mathop{\rm Var}\nolimits}

\newcommand{\pff}{\noindent {\sc Proof.}\ }
\newcommand{\fin}{\sqcap\!\!\!\!\sqcup}

\newcommand\1{\leavevmode\hbox{\rm \small1\kern-0.35em\normalsize1}}

\def\g#1{\mathbb #1}
\def\rond#1{\mathcal #1}

\def\C{\g C}

\newtheorem {Th}{Theorem}
\newtheorem {Cor}{Corollary}

\newtheorem {Lem}{Lemma}
\newtheorem {Rem}{Remark}

\begin{document}

\begin{center}
{\bf 
B-urns\footnote{

Fringe analysis.
P\'olya urn.
Urn model.
Martingale.
Multitype branching process.
Smoothing transforms.
Contraction method.
}

}
\end{center}

\begin{center}
{\fontauthors
Brigitte Chauvin $ ^a$,
%\footnote{INRIA Rocquencourt, project Algorithms and Laboratoire de Math\'ematiques de Versailles, CNRS, UMR 8100 --  INRIA Domaine de Voluceau B.P.105, 78153 Le Chesnay CEDEX (France).},
Dani\`ele Gardy $ ^b$, 
Nicolas Pouyanne $ ^a$ and
Dai-Hai Ton-That $ ^b$
%\footnote{Laboratoire de Math\'ematiques de Versailles, CNRS, UMR 8100 -- Universit\'e de Versailles - St-Quentin,  45 avenue des Etats-Unis, 78035 Versailles CEDEX (France).
}

\medskip

$ ^a$ Laboratoire de Math\'ematiques de Versailles, CNRS UMR 8100. \\
$ ^b$ Laboratoire PRiSM, CNRS UMR 8144.\\
 Universit\'e de Versailles - St-Quentin, \\
45, avenue des Etats-Unis, 78035 Versailles Cedex, France.
\end{center}

\vskip 20pt
\begin{center}
July 22nd, 2015\\
%by
\end{center}

\vskip 1truecm
\noindent{\bf Abstract.}

The fringe of a B-tree with parameter $m$ is considered as a particular P\'olya urn with $m$ colors. More precisely, the asymptotic behaviour of this fringe, when the number of stored keys tends to infinity, is studied through the composition vector of the fringe nodes. We establish its typical behaviour together with the fluctuations around it. The well known phase transition in P\'olya urns has the following effect on  B-trees: for $m\leq 59$, the fluctuations are asymptotically Gaussian, though for $m\geq 60$, the composition vector is oscillating;  after scaling, the fluctuations of such an urn strongly converge to a random variable $W$. This limit is $\g C$-valued and it does not seem to follow any classical law. Several properties of $W$ are shown: existence of exponential moments, characterization of its distribution as the solution of a smoothing equation, existence of a density relatively to the Lebesgue measure on $\g C$, support of $W$. Moreover, a few representations of the composition vector for various values of $m$ illustrate the different kinds of convergence.

\tableofcontents
%%%%%%%%%%%%%%%%%%%%%%%%%%
\section{Introduction}
\label{intro}

B-trees are a fundamental structure in computer science, they have been introduced in the early seventies by Bayer and McCreight \cite{Bay, BayMcC},
to store large quantities of data.  These particular search trees are conceived in order to have all their leaves at the same level. The nodes at the deepest level are called the \emph{fringe nodes}. A precise description can be found in Section \ref{sec-algo} where are presented two classical algorithms giving a B-tree. The actual writing of these algorithms can be found for example in Cormen et al. \cite{CorLeiRiv} for one of them (the so-called \emph{prudent} algorithm in the sequel), in Kruse and Ryba \cite{KruRyb} for the other one (called the \emph{optimistic} algorithm in the sequel).

\medskip

The fringe analysis of B-trees goes back to Yao \cite{Yao} and has been developped by many authors (see for example the Baeza-Yates' survey \cite{BaeSurvey}), both for B-trees and B$^{+}$-trees (where all the keys are stored in the fringe nodes). In Yao's paper \cite{Yao} appears the P\'olya urn model, which we develop in this article.
Indeed, the fringe of a B-tree with parameter $m$ (where $m$ is a positive integer) can be considered as a particular P\'olya urn with $m$ colors, so that a lot of information can be obtained concerning the asymptotic behaviour of this fringe, when the number of stored keys tends to infinity.

 \medskip
 
Let us describe a P\'olya urn process as follows. 
Consider an urn that contains balls of, say, $d$ different colors. 
Start with a finite number of different color balls as initial composition (possibly monochromatic).
At each discrete time $n$, draw a ball at random, check its color,
put it back into the urn and add balls according to the following rule:
if the drawn ball  is of color $i$,  add $a_{i,j}$  balls of color $j$, where the $a_{i,j}$ are integer-valued. 
Thus, the replacement rule is described by the so-called \emph{replacement matrix}, which is a dimension $d$ matrix, whose coefficients are the  $a_{i,j}$, for $i$ and $j$ in $\{ 1, \dots , d\}$.

Usually, the integers $a_{i,j}$ are assumed to be nonnegative for $i\not= j$ and the integers $a_{i,i}$ are nonpositive or nonnegative. A negative coefficient $a_{i,i}$ means that, if a ball of color $i$ is drawn, then $a_{i,i}$ balls of color $i$ are removed from the urn. In this case, we have to ensure that at least $a_{i,i}$ balls of color $i$ exist in the urn. This quality is called the \emph{tenability} of the urn.
 %  TENABILITY
 % ne pas citer le livre de Mahmoud, page 49, c'est trop mal fait : ˆ deux couleurs seulement, et sans dire que les conditions arithmŽtiques sont nŽcessaires et suffisantes.
To ensure that an urn with a negative coefficient $a_{i,i}$ is tenable, it is necessary and sufficient to have the following arithmetical condition (this can be easily proved by  induction on $n$). Fix an initial composition $(\alpha_1, \dots , \alpha_d)$,
meaning that there are $\alpha_j$ balls of colour $j$ at time zero in the urn, then the tenability condition can be written as
\begin{equation}
\label{tenable}
- a_{i,i} \hbox{ divides }  \alpha_i, a_{1,i}, \dots , a_{d,i}  .
\end{equation}
Moreover, in the present paper, the urn is assumed to be {\it balanced}, which means that the total number of balls added at each
step is a constant: there exists an integer $S$ such that, for any $i$ in $\{ 1, \dots , d\}$, $\displaystyle\sum_{j=1}^d a_{i,j} = S$.

Let us emphasize that ``drawing a ball at random'' means choosing \emph{uniformly} among the balls contained in the urn.
That is why this model is related to many situations in mathematics, algorithmics or theoretical physics where a uniform choice among objects determines the evolution of a process. See Johnson and Kotz's book \cite{JK}, Mahmoud's book \cite{Mah08} or Flajolet et al.  \cite{FlaDumPuy} for many examples. 
For a general probabilistic treatment of P\'olya urns, see \cite{Pou08}, Janson \cite{Jan} or Mailler \cite{Mailler}.
 
   \medskip
  
 In Yao's paper \cite{Yao}, the focus is on the average number of nodes in the B-tree. Nevertheless, the main ideas are already there, namely the dynamics transforming a tree of size $n$ into a tree of size $n+1$, which is the same dynamics as in a  P\'olya urn process. The recent progresses in P\'olya urn processes and their asymptotic behaviour (\cite{Pou08, ChaPouSah, ChaMaiPou}, Janson \cite{Jan},  Mailler \cite{Mailler}) lead to a more complete landscape for the B-trees. Our aim in this article is to present in a hopefully concise form a collection of results about the asymptotic behaviour of the fringe nodes in a B-tree, namely their typical behaviour and the fluctuations around it. Our main interest is focused on these fluctuations, which happen to have a  phase transition: for $m\leq 59$, the fluctuations are of order $\sqrt n$ and have a Gaussian limit in distribution. But for $m\geq 60$, the fluctuations  are of order $n^{\sigma}$, where $\sigma$ is larger than $1/2$ and increases to~$1$ when $m$ tends to infinity. Moreover, an oscillating and significative phenomenon occurs in the fluctuation term.
% are of order $n^{\lambda}$ where $\lambda$ is a complex, nonreal number, whose real part increases to $1$ when $m$ increases. Consequently, due to the imaginary part of $\lambda$,  an oscillating, significative behaviour occurs. 
After scaling, the fluctuations strongly converge (meaning almost surely) to a random limit, here called $W$. The random variable $W$ is $\g C$-valued and does not seem to follow any classical law. 
 
\medskip

The paper is organized as follows. In Section \ref{sec-algo} are presented  two classical algorithms allowing to construct a B-tree. In that section is also precised how the insertion dynamics is that of a suitable P\'olya urn. In Section \ref{sec-ProcessusDiscrets} are introduced the random vectors which describe the fringe of a B-tree. In Section \ref{sec-transition} is established the phase transition, and we get the precise asymptotic behaviour of the fringe nodes, in Corollary~\ref{asymptoticsBtreeL} of Theorem \ref{asymptoticsDTL}. For $m\geq 60$, the fluctuations around the drift are expressed via a random variable $W$, which is studied in the last sections. 
Thanks to an embedding into continuous time (Section \ref{sec-embedding}), a multitype branching process is put forward. Properties of the continuous-time limit process can be translated to the discrete-time process, via an explicit connection. 
Several properties of $W$ are proved in Section \ref{sec-limitlaw}: $W$ admits a density on the whole complex plane; it  has exponential moments; it is the unique solution of a certain ``smoothing equation''  in a convenient probability distribution space. 
Finally, in Section \ref{sec-simulations}, a few pictures provide a synthetic and concrete illustration of the different kinds of convergence, depending on whether $m\leq 59$ or $m\geq 60$.
%Nevertheless, some unanswered questions remain, gathered in Section \ref{sec-conclusion}.

%%%%%%%%%
%%  pour la version "journal", il serait mieux d'Žécrire une intro plus classique pour les matheux, ˆ savoir mettre les réŽsultats pour Žéviter de lire le papier.

%%%%%%%%%%%%%%%%%%%%%%%%%
\section{B-tree algorithms}
\label{sec-algo}

\subsection{Description of a B-tree}

For a positive integer $m\geq 2$, a B-tree with parameter $m$ is a search\footnote
{A search tree is a tree where internal nodes contain sorted keys and where a node containing $k$ keys $x_1, \dots , x_k$ defines $k+1$ intervals such that, for $j=1, \dots , k+1$, the keys in the $j$-th subtree belong to the $j$-th interval.} 
tree, where the keys are stored into the internal nodes and the leaves\footnote
{The leaves, sometimes called external nodes, are the nodes without any descendant.} 
represent insertion possibilities (we call them \emph{gaps}), they do not contain any key;
furthermore all the leaves are at the same depth. 
A \emph{fringe node} is an internal node whose only descendants are  leaves. In the literature, these fringe nodes are sometimes called \emph{final internal nodes} or \emph{leaf-nodes}, or \emph{internal leaves}. We try to be non-ambiguous in the following, and use the terms fringe nodes and fringe node process. In the figures below, internal nodes are represented by ellipses and leaves by squares.

 As is the case for the leaves, the fringe nodes of a B-tree are at the same depth.
Moreover, each internal node (fringe or otherwise) has a capacity; the root contains between $1$ and  $C(m)$ keys, and the other internal nodes between $c(m)$ and $C(m)$ keys. 
When a node contains $C(m)$ keys, we say that the node is \emph{saturated}.

\medskip
 
The minimal -- $c(m)$ -- and maximal -- $C(m)$ -- values depend both on the parameter $m$ and on the precise definition of the B-tree, which is itself closely related to the exact insertion algorithm, of which we present two versions below.
Let us just state that $c(m) = m-1$ and $C(m) = 2m-1$ for the first algorithm, and  $c(m) = m$ and $C(m) = 2m$ for the second one.
In both cases we want to insert a new key into a tree of size $n$, i.e. having already $n$ keys in its internal nodes, and consequently $n+1$ leaves, or insertion possibilities.

\vskip 5mm

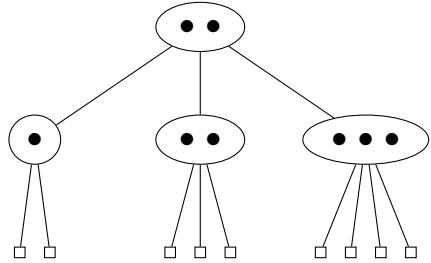
\begin{figure}[h!]
\begin{center}
%%%%%% arbre B basique pour m=3 algo prudent
\tikzstyle{interne}=[ellipse,draw,dashed, fill=gray!20]
\tikzstyle{entree}=[ellipse,draw]
\tikzstyle{entreerose}=[ellipse,fill=pink,draw=pink]
\tikzstyle{entreerouge}=[ellipse,fill=red,draw=red]
\tikzstyle{gap}=[rectangle,draw,inner sep=2pt]
\tikzstyle{gaprose}=[rectangle,draw,fill=pink,inner sep=2pt]
\tikzstyle{gapvert}=[rectangle,draw,fill=green,inner sep=2pt]
\tikzstyle{gaprouge}=[rectangle,draw,fill=red,inner sep=2pt]
\tikzstyle{pointeur}=[>=latex,<-,thick,level distance = 8mm]
\begin{tikzpicture}
[level 1/.style={sibling distance=22mm},
level 2/.style={sibling distance=4mm},
level 3/.style={sibling distance=4mm},
noedge/.style={edge from parent/.style={}},
    nonode/.style={}
]
\node [entree]  at (3,0) {$\bullet$ $\bullet$}
child {node [entree]  { $\bullet$}
	child{node [gap] {}}
	child{node [gap] {}}
  }
  child {node [entree]  {$\bullet$ $\bullet$}
  	child{node [gap] {}}
	child{node [gap] {}}
	child{node [gap] {}}
  }
  child {node [entree]  {$\bullet$ $\bullet$ $\bullet$}
  	child{node [gap] {}
		}
	child{node [gap] {}}
	child{node [gap] {}}
	child{node [gap] {}}
  }
 ;
\end{tikzpicture}
%\vskip -20pt
\caption{%
\label{fig:Btree-basic} {\it An example of B-tree of size $8$. Here $m=2$, nodes contain between $1$ and $3$ keys. There are $3$ fringe nodes and $9$ leaves. One fringe node is saturated.
}}
\end{center}
\end{figure}

\subsection{The prudent algorithm}

\begin{figure}[h!]
%%%%%% arbre B sans couleurs pour m=3 algo standard 
\tikzstyle{interne}=[ellipse,draw,dashed, fill=gray!20]
\tikzstyle{entree}=[ellipse,draw]
\tikzstyle{entreerose}=[ellipse,fill=pink,draw=pink]
\tikzstyle{entreerouge}=[ellipse,fill=red,draw=red]
\tikzstyle{gap}=[rectangle,draw,inner sep=2pt]
\tikzstyle{gaprose}=[rectangle,draw,fill=pink,inner sep=2pt]
\tikzstyle{gapvert}=[rectangle,draw,fill=green,inner sep=2pt]
\tikzstyle{gaprouge}=[rectangle,draw,fill=red,inner sep=2pt]
\tikzstyle{pointeur}=[>=latex,<-,thick,level distance = 8mm]
\begin{tikzpicture}
[level 1/.style={sibling distance=22mm},
level 2/.style={sibling distance=4mm},
level 3/.style={sibling distance=4mm},
noedge/.style={edge from parent/.style={}},
    nonode/.style={}
]
\node [entree]  at (3,0) {$\bullet$ $\bullet$}
child {node [entree]  { $\bullet$}
	child{node [gap] {}}
	child{node [gap] {}}
  }
  child {node [entree]  {$\bullet$ $\bullet$}
  	child{node [gap] {}}
	child{node [gap] {}}
	child{node [gap] {}}
  }
  child {node [entree]  {$\bullet$ ${\color{red}\bullet}$ $\bullet$}
  	child{node [gapvert] {}
		child[pointeur]{node [nonode]{}}
		child[noedge]{node [nonode]{}}
		}
	child{node [gap] {}}
	child{node [gap] {}}
	child{node [gap] {}}
  }
 ;
 \node [entree]  at (12,0) {$\bullet$ $\bullet$ ${\color{red}\bullet}$}
child {node [entree]  {$\bullet$}
  	child{node [gap] {}}
	child{node [gap] {}}
  }
  child {node [entree]  {$\bullet$ $\bullet$}
  	child{node [gap] {}}
	child{node [gap] {}}
	child{node [gap] {}}
  }
  child {node [entree]  {${\color{green}\bullet}$ $\bullet$}
  	child{node [gap] {}}
	child{node [gap] {}}
	child{node [gap] {}}
  }
  child {node [entree]  {$\bullet$}
  	child{node [gap] {}}
	child{node [gap] {}}
  }
 ;
\end{tikzpicture}
\vskip -20pt
\caption{%
\label{fig:Btree-standardNC} {\it An example of insertion in a B-tree, for the prudent algorithm. Here $m=2$, nodes contain between $1$ and $3$ keys. The middle key in red moves to the parent node.
}}
\end{figure}
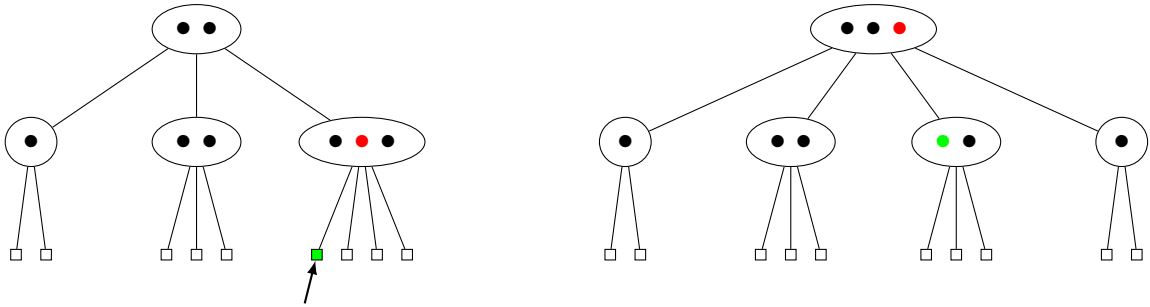

In what we call here the \emph{prudent algorithm} for insertion into a B-tree with parameter $m$, the nodes contain between $m-1$ and $2m-1$ keys. An insertion of a new key concerns a given leaf, so that a branch (the nodes between the root and this leaf) is determined for this insertion. The algorithm proceeds by \emph{going down} from the root to the leaf, along this branch.
We begin by checking the root: if it is saturated, it is split, a new root is created with a single key which is the median of the keys  of the old root (remember that it has an odd number of keys, hence the median is defined without any ambiguity) and two sons, and the height of the tree increases by $1$. If the root is not saturated, we do not modify it.
We then proceed along the branch to the insertion gap.
When we meet a  (non-root) saturated node, the median key of that node  moves to the parent node (which is not saturated -- if it  initially was, we have already taken care of it) and the saturated node is split. 
Then, when we finally arrive at a fringe node, we split it when necessary, and the insertion of the new key always takes place into a non-saturated fringe node: the saturated nodes are dealt with \emph{before} we find the node in which the insertion of the new key will take place.
This algorithm, which can be presented both recursively and iteratively (there being only a descent from the root to a leaf), is found, e.g., in the book of Cormen et al \cite{CorLeiRiv}. 
If we consider the fringe nodes, insertion on a saturated node (with $2m-1$ keys) gives rise to $2$ new fringe nodes with respectively $m$ and $m-1$ keys. See Figure \ref{fig:Btree-standardNC}.

\vskip 5mm

\subsection{The optimistic algorithm}

\begin{figure}[h]
%%%%%% arbre B sans couleurs pour m=3 algo grandiose
\tikzstyle{interne}=[ellipse,draw,dashed, fill=gray!20]
\tikzstyle{entree}=[ellipse,draw]
\tikzstyle{entreerose}=[ellipse,fill=pink,draw=pink]
\tikzstyle{entreerouge}=[ellipse,fill=red,draw=red]
\tikzstyle{gap}=[rectangle,draw,inner sep=2pt]
\tikzstyle{gaprose}=[rectangle,draw,fill=pink,inner sep=2pt]
\tikzstyle{gapvert}=[rectangle,draw,fill=green,inner sep=2pt]
\tikzstyle{gaprouge}=[rectangle,draw,fill=red,inner sep=2pt]
\tikzstyle{pointeur}=[>=latex,<-,thick,level distance = 8mm]
\begin{tikzpicture}
[level 1/.style={sibling distance=22mm},
level 2/.style={sibling distance=4mm},
level 3/.style={sibling distance=4mm},
noedge/.style={edge from parent/.style={}},
    nonode/.style={}
]
\node [entree]  at (3,0) {$\bullet$ $\bullet$}
child {node [entree]  {$\bullet$ $\bullet$}
  	child{node [gap] {}}
	child{node [gap] {}}
	child{node [gap] {}}
  }
  child {node [entree]  {$\bullet$ $\bullet$ $\bullet$}
  	child{node [gap] {}}
	child{node [gap] {}}
	child{node [gap] {}}
	child{node [gap] {}}
  }
  child {node [entree]  {$\bullet$ ${\color{red}\bullet}$ $\bullet$ $\bullet$}
  	child{node [gapvert] {}
		child[pointeur]{node [nonode]{}}
		child[noedge]{node [nonode]{}}
		}
	child{node [gap] {}}
	child{node [gap] {}}
	child{node [gap] {}}
	child{node [gap] {}}
  }
 ;
 \node [entree]  at (12,0) {$\bullet$ $\bullet$ ${\color{red}\bullet}$}
child {node [entree]  {$\bullet$ $\bullet$}
  	child{node [gap] {}}
	child{node [gap] {}}
	child{node [gap] {}}
  }
  child {node [entree]  {$\bullet$ $\bullet$ $\bullet$}
  	child{node [gap] {}}
	child{node [gap] {}}
	child{node [gap] {}}
	child{node [gap] {}}
  }
  child {node [entree]  {${\color{green}\bullet}$ $\bullet$}
  	child{node [gap] {}}
	child{node [gap] {}}
	child{node [gap] {}}
}
child {node [entree]  {$\bullet$ $\bullet$}
  	child{node [gap] {}}
	child{node [gap] {}}
	child{node [gap] {}}
  }
 ;
\end{tikzpicture}
\vskip -20pt
\caption{%
\label{fig:BtreeNC} {\it An example of insertion in a B-tree, for the optimistic algorithm. Here $m=3$, nodes contain between $2$ and $4$ keys. The middle key in red is determined among the $4$ keys of the saturated node and the new key, and moves to the parent node.
}}
\end{figure}
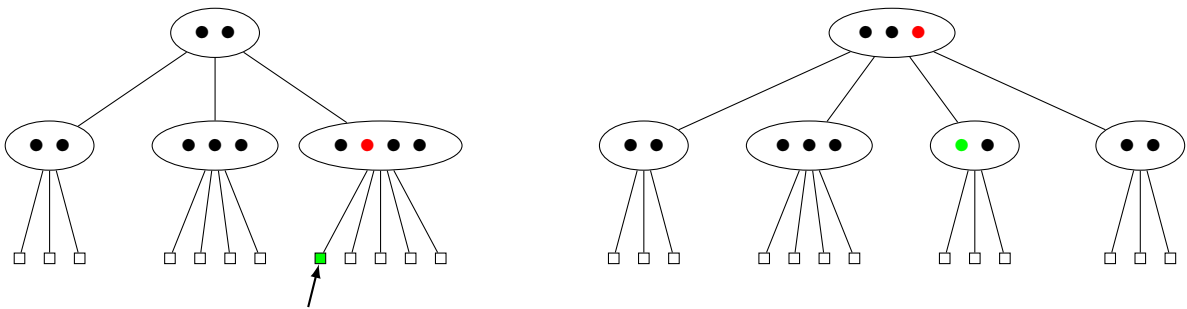

In what we call here the \emph{optimistic algorithm}, for insertion into a B-tree with parameter $m-1$, the nodes contain between $m-1$ and $2m-2$ keys. 
Here the saturated nodes are dealt with \emph{after}  we have found the place for insertion.
An insertion of a new key concerns a given leaf. 
If the corresponding fringe node is not saturated, the insertion occurs in this node; if it is saturated, the algorithm has to create a non-saturated node into which we can insert the new key.
It needs to find what would be the place of the new key among the (already sorted) $2m-2$ keys; the middle key among these $2m-2 +1 = 2m-1$ keys moves to the parent node, and the saturated node is split into $2$ new fringe nodes with $m-1$ keys. If the parent node is saturated, a key is pushed up into the grandparent node, etc... all the way up to the root if necessary; if the root is saturated, it is split as well and the height of the tree increases by $1$. This algorithm proceeds by going down from the root to the gap of insertion, and then up to (some node on) the branch from that leaf to the root, and is possibly best understood recursively. Figure \ref{fig:BtreeNC} illustrates an insertion on a saturated node for a B-tree with parameter $m=3$.

\vskip 5mm

\subsection{Insertion as the evolution of a P\'olya urn}

For both the prudent and the optimistic algorithms, let us define different types of fringe nodes: we say that a fringe node is \emph{of type $k$}, when it contains $m+k-2$ keys and has thus $m+k-1$ gaps. For the prudent algorithm, $k$  varies between $1$ and $m+1$;  for the optimistic algorithm, $k$  varies between $1$ and $m$ . 

\medskip

We analyse the fringe of the tree through the so-called \emph{composition vector} $L_n$, which counts the number of fringe nodes of each type in a B-tree of parameter $m$, at time~$n$, i.e., assuming that we start from an empty tree and add the keys one by one, when the tree contains $n$ keys. Thus, $L_n^{(k)}$, the $k$-th coordinate of $L_n$, is the number of fringe nodes of type $k$. For the prudent algorithm, $L_n$ is a vector of dimension $m+1$, whereas it is a vector of dimension $m$ for the optimistic algorithm .

\medskip

For both algorithms, we  define $G_n$ as the composition vector of \emph{gaps} at time $n$. We say that a gap is \emph{of type $k$}, when it is attached to a fringe node of type $k$. Thus $G_n^{(k)}$, the $k$-th coordinate of $G_n$, is the number of gaps of type $k$. In other words: 
\begin{equation}
\label{gap-node}
(m+k-1) L_n^{(k)} = G_n^{(k)}.
\end{equation}

\medskip

For both algorithms, the process $(G_n)_{n\in\g N}$ is a P\'olya urn process, as defined in the Introduction, where the balls are the gaps and the colors are the different types. Indeed, when the keys are randomly chosen under the so-called random permutation model, meaning that the keys are independently identically distributed (i.i.d.), then the insertion of a new key in a B-tree of size $n$ occurs \emph{uniformly} on any of the $n+1$ gaps of the tree. 

\medskip

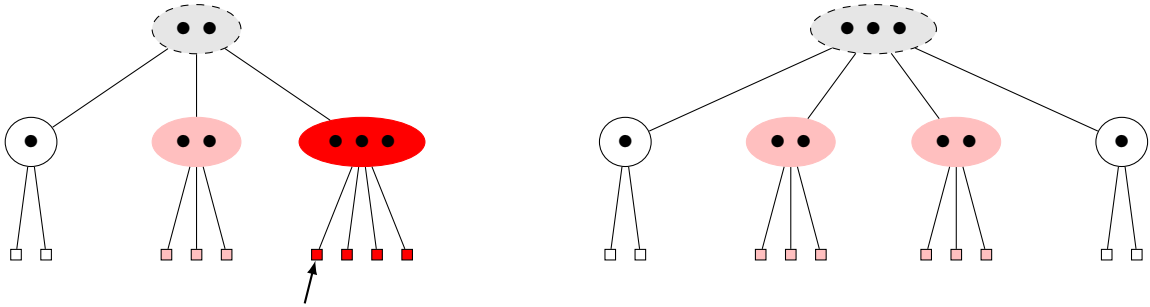
\begin{figure}[h]
%%%%%% arbre B avec couleurs pour m=3 algo standard 
%\begin{center}
\tikzstyle{interne}=[ellipse,draw,dashed, fill=gray!20]
\tikzstyle{entree}=[ellipse,draw]
\tikzstyle{entreerose}=[ellipse,fill=pink,draw=pink]
\tikzstyle{entreerouge}=[ellipse,fill=red,draw=red]
\tikzstyle{gap}=[rectangle,draw,inner sep=2pt]
\tikzstyle{gaprose}=[rectangle,draw,fill=pink,inner sep=2pt]
\tikzstyle{gaprouge}=[rectangle,draw,fill=red,inner sep=2pt]
\tikzstyle{pointeur}=[>=latex,<-,thick,level distance = 8mm]
\begin{tikzpicture}
[level 1/.style={sibling distance=22mm},
level 2/.style={sibling distance=4mm},
level 3/.style={sibling distance=4mm},
noedge/.style={edge from parent/.style={}},
    nonode/.style={}
]
\node [interne]  at (3,0) {$\bullet$ $\bullet$}
child {node [entree]  { $\bullet$}
	child{node [gap] {}}
	child{node [gap] {}}
  }
  child {node [entreerose]  {$\bullet$ $\bullet$}
  	child{node [gaprose] {}}
	child{node [gaprose] {}}
	child{node [gaprose] {}}
  }
  child {node [entreerouge]  {$\bullet$ $\bullet$ $\bullet$}
  	child{node [gaprouge] {}
		child[pointeur]{node [nonode]{}}
		child[noedge]{node [nonode]{}}
		}
	child{node [gaprouge] {}}
	child{node [gaprouge] {}}
	child{node [gaprouge] {}}
  }
 ;
 \node [interne]  at (12,0) {$\bullet$ $\bullet$ $\bullet$}
child {node [entree]  {$\bullet$}
  	child{node [gap] {}}
	child{node [gap] {}}
  }
  child {node [entreerose]  {$\bullet$ $\bullet$}
  	child{node [gaprose] {}}
	child{node [gaprose] {}}
	child{node [gaprose] {}}
  }
  child {node [entreerose]  {$\bullet$ $\bullet$}
  	child{node [gaprose] {}}
	child{node [gaprose] {}}
	child{node [gaprose] {}}
  }
  child {node [entree]  {$\bullet$}
  	child{node [gap] {}}
	child{node [gap] {}}
  }
 ;
\end{tikzpicture}
%\end{center}
\vskip -20pt
\caption{%
\label{fig:Btree-standard} {\it An example of insertion in a B-tree, for the prudent algorithm. Here $m=2$, nodes contain between $1$ and $3$ keys, there are $3$ colors, white, pink and red.
}}
\end{figure}

$\bullet$ In the prudent algorithm, the number of keys in a fringe node ranges from $m-1$ to $2m-1$, there are $m+1$ types, and the vector $L_n$ is of dimension $m+1$. The replacement matrix of the gap process is of dimension $m+1$ and equal to\footnote{An empty entry stands for a zero in all the matrices of this article.}
$$
r_m=
\left(
\begin{array}{ccccc}
-m&(m+1)&&&\\
&-(m+1)&(m+2)&&\\
&&\ddots&\ddots&\\
%&&\ddots&\ddots&\\
&&&\ddots&2m\\
m&(m+1)&&&-2m
\end{array}
\right).
$$
 Figure \ref{fig:Btree-standard} illustrates the same insertion as in Figure~\ref{fig:Btree-standardNC}, taking into account the different types (colors) of the fringe nodes.

\medskip

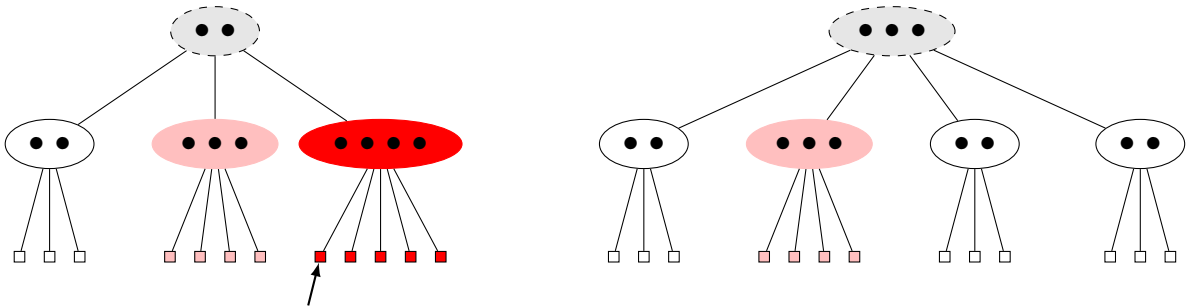
\begin{figure}[h]
\begin{center}
%%%%%% arbre B avec couleurs pour m=3 algo grandiose
\tikzstyle{interne}=[ellipse,draw,dashed, fill=gray!20]
\tikzstyle{entree}=[ellipse,draw]
\tikzstyle{entreerose}=[ellipse,fill=pink,draw=pink]
\tikzstyle{entreerouge}=[ellipse,fill=red,draw=red]
\tikzstyle{gap}=[rectangle,draw,inner sep=2pt]
\tikzstyle{gaprose}=[rectangle,draw,fill=pink,inner sep=2pt]
\tikzstyle{gaprouge}=[rectangle,draw,fill=red,inner sep=2pt]
\tikzstyle{pointeur}=[>=latex,<-,thick,level distance = 8mm]
\begin{tikzpicture}
[ level 1/.style={sibling distance=22mm},
level 2/.style={sibling distance=4mm},
level 3/.style={sibling distance=4mm},
noedge/.style={edge from parent/.style={}},
    nonode/.style={}
]
\node [interne]  at (3,0) {$\bullet$ $\bullet$}
child {node [entree]  {$\bullet$ $\bullet$}
  	child{node [gap] {}}
	child{node [gap] {}}
	child{node [gap] {}}
  }
  child {node [entreerose]  {$\bullet$ $\bullet$ $\bullet$}
  	child{node [gaprose] {}}
	child{node [gaprose] {}}
	child{node [gaprose] {}}
	child{node [gaprose] {}}
  }
  child {node [entreerouge]  {$\bullet$ $\bullet$ $\bullet$ $\bullet$}
  	child{node [gaprouge] {}
		child[pointeur]{node [nonode]{}}
		child[noedge]{node [nonode]{}}
		}
	child{node [gaprouge] {}}
	child{node [gaprouge] {}}
	child{node [gaprouge] {}}
	child{node [gaprouge] {}}
  }
 ;
 \node [interne]  at (12,0) {$\bullet$ $\bullet$ $\bullet$}
child {node [entree]  {$\bullet$ $\bullet$}
  	child{node [gap] {}}
	child{node [gap] {}}
	child{node [gap] {}}
  }
  child {node [entreerose]  {$\bullet$ $\bullet$ $\bullet$}
  	child{node [gaprose] {}}
	child{node [gaprose] {}}
	child{node [gaprose] {}}
	child{node [gaprose] {}}
  }
  child {node [entree]  {$\bullet$ $\bullet$}
  	child{node [gap] {}}
	child{node [gap] {}}
	child{node [gap] {}}
}
child {node [entree]  {$\bullet$ $\bullet$}
  	child{node [gap] {}}
	child{node [gap] {}}
	child{node [gap] {}}
  }
 ;
\end{tikzpicture}
\end{center}
\vskip -30pt
\caption{%
\label{fig:Btree} {\it An example of insertion in a B-tree, for the optimistic algorithm. Here $m=3$, nodes contain between $2$ and $4$ keys, there are $3$ colors, white, pink and red.
}}
\end{figure}

$\bullet$ In the optimistic algorithm, the number of keys in a fringe node ranges from $m-1$ to $2m-2$, there are $m$ types, and the vector $L_n$ is of dimension $m$. The replacement matrix of the gap process is of dimension $m$ and equal to
$$
R_m=\begin{pmatrix}
-m&m+1&&&\\
&-(m+1)&m+2&&\\
&&\ddots&&\\
&&&-(2m-2)&2m-1\\
2m&&&&-(2m-1)
\end{pmatrix}.
$$
Figure \ref{fig:Btree} illustrates the same insertion as in Figure \ref{fig:BtreeNC}, taking into account the different types (colors) of the fringe nodes.

\medskip

Observe that both replacement matrices  $r_m$ and $R_m$ are balanced (any row sums to~$1$), which is an immediate consequence of the dynamics, since one key (one ball) is added at each unit of time. 

\medskip

All the results in this paper hold for both algorithms, including the phase transition depending on whether $m\leq 59$ or $m\geq 60$.
However, the proofs and results for the optimistic algorithm being somewhat simpler than those for the prudent algorithm, we choose to present them and to leave the other case to the reader: from now on, 
\begin{center}
\emph{we consider a B-tree constructed by the optimistic algorithm.}
\end{center}

%%%%%%%%%%%%%%%%%%%%%%%%%
\section{Gaps of a B-tree as a P\'olya urn}
\label{sec-ProcessusDiscrets}

Let us remind from Section \ref{sec-algo} that a fringe node of a B-tree contains from $m-1$ to $2m-2$ keys and from $m$ to $2m-1$ gaps.
For any $k\in\{ 1,\dots ,m\}$, a fringe node that contains $m+k-2$ keys is called \emph{of type $k$}. We are interested in the \emph{fringe node composition vector} $L_n$ of a B-tree at time $n$, whose $k$-th coordinate counts the number of fringe nodes of type $k$.

\vskip 20pt
{\bf Notation}

Let $\left( e_k\right) _{1\leq k\leq m}$ be the canonical basis of $\g R^m$.
Denote by $w_1,\dots ,w_m$ the vectors defined by
\begin{equation}
\label{increments}
\left\{
\begin{array}{l}
w_k=-e_k+e_{k+1},~~1\leq k\leq m-1\\
w_m=2e_1-e_{m}.
\end{array}
\right.
\end{equation}

The $w_k$'s are the increment vectors of the fringe node dynamics:
when a key is inserted in a fringe node of type $k\in\{ 1,\dots ,m-1\}$, the fringe node is replaced by a fringe node of type $k+1$ (addition of vector $w_k$) and when a key is inserted in a fringe node of type $m$, the fringe node is replaced by two fringe nodes of type $1$ (addition of vector $w_m$).
When the keys are randomly drawn under the permutation model, the insertion is \emph{uniform} on the gaps and the fringe node composition process of a B-tree is modelized by the $\g R^m$-valued Markov chain $\left( L_n\right) _{n\in\g N}$ defined as follows by its transition probabilities.

\vskip 20pt
{\bf Definition of the (discrete-time) fringe node process}

For any $k\in\{ 1,\dots ,m\}$,
\begin{equation}
\label{markovLeaves}
\Proba\left( L_{n+1}=L_n+w_k\Big| L_n\right)
=\frac{(m+k-1)\langle L_n,e_k\rangle}{K_n}
\end{equation}
where the scalar product $\langle L_n,e_k\rangle = L_n^{(k)}$ is the $k$-th coordinate of $L_n$ and where $K_n$ denotes the total number of gaps at time $n$.
Of course, when the process starts initially with $N_0$ keys at time $0$, then $K_n=1+N_0+n$.
Note that, considering the B-tree, the number of gaps in a fringe node of type $k$ is $m+k-1$ whereas the total number of gaps in the tree at time $n$ is exactly $K_n$ so that Formulae~\eqref{markovLeaves} completely define a Markov process;
it reflects the uniform insertions of the keys in the gaps.

\vskip 10pt
Alternatively, one can consider the gap process.
A gap is called of type $k$ when it is contained in a fringe node of type $k$.
Note once more that a fringe node of type $k$ contains $m+k-1$ gaps.
Expressed in terms of gaps, the dynamics of key insertion in the B-tree is the following:
when the key is inserted in a gap of type $k\in\{ 1,\dots ,m-1\}$, then $m+k-1$ gaps of type $k$ disappear and are replaced by $m+k$ gaps of type $k+1$;
when the key is inserted in a gap of type $m$, then $2m-1$ gaps of type $m$ disappear and are replaced by $2m$ gaps of type $1$.
Moreover, under the random permutation model, \emph{the gaps are drawn uniformly}.
In other words, the gap composition process of a B-tree is modelized by the following $m$-color P\'olya urn process $\left( G_n\right) _{n\in\g N}$, having $R_m$ as replacement matrix and $(m,0,\dots ,0)$ as initial composition.

%V2 This shows that, considering their different types, the gap process is exactly the composition process of the following P\'olya urn denoted by $\left( G_n\right) _{n\in\g N}$. In the following, we consider on the same probability space the physical gap process and the P\'olya urn process and, abusing notations, we denote both processes by the same notation $\left( G_n\right) _{n\in\g N}$.

\vskip 20pt
{\bf Definition of the (discrete-time) gap process}

Denote by $\left( G_n\right) _{n\in\g N}$ the $m$-color P\'olya urn process defined by the $m$-dimensional replacement matrix 
$$
R_m=\begin{pmatrix}
-m&m+1&&&\\
&-(m+1)&m+2&&\\
&&\ddots&&\\
&&&-(2m-2)&2m-1\\
2m&&&&-(2m-1)
\end{pmatrix}.
$$
Its balance, namely its common row sum, equals $1$.
Note that the diagonal entries are negative.
Nevertheless, the urn is tenable because, in any column, all entries are multiple of the diagonal coefficient:
when a ball of color $2$ is drawn, $m+1$ extra balls must be withdrawn from the urn which is always possible because balls of type $2$ are put $m+1$ by $m+1$ in the urn as can be seen on $R_m$'s second column.
The same phenomenon occurs for any color.
Of course, these negative diagonal entries  imply that one must necessarily take an initial composition that satisfies such divisibility conditions as well:
for any $k\in\{ 1,\dots ,m\}$, the initial number of balls of type $k$ satisfies
$$
m+k-1~\hbox{ divides }~\langle L_0,e_k\rangle.
$$
The symbol $\langle .,. \rangle$ denotes here the standard scalar product on $\g R^m$.
Thus, the condition \eqref{tenable} is fulfilled.

\vskip 10pt
Both Markov processes $\left( L_n\right) _{n}$ and $\left( G_n\right) _{n}$ are related by Relation~\eqref{lienLeafGapDiscret}, stated hereunder.
Let $P$ be the $m$-dimensional diagonal invertible matrix
\begin{equation}
\label{defP}
P=\Diag\left( m,m+1,\dots ,2m-1\right).
\end{equation}
When $V\in\g N^m\setminus\{ 0\}$, denote by $\left( L_n^V\right) _{n\geq 0}$ the fringe node process starting with $L_0=V$ and by $\left( G_n^V\right) _{n\geq 0}$ the gap process starting from $G_0=V$.
Then, for any $V\in\g N^m\setminus\{ 0\}$, one has\footnote{
When no confusion is possible, we denote by PV the product of the square matrix $P$ by the vector $V\in\g R^m$ instead of the correct form $P\transp V$.
}
immediatly from \eqref{gap-node}
\begin{equation}
\label{lienLeafGapDiscret}
\left( G_n^{PV}\right) _{n\in\g N}\egalLoi \left( PL_n^V\right) _{n\in\g N}.
\end{equation}
In particular, denoting by $|V|$ the sum of $V$'s coordinates, the total number of gaps in the B-tree at time $n$ is
$$
K_n=n+|PV|
=\sum _{k=1}^m(m+k-1)\langle L_n^{V},e_k\rangle
=\sum _{k=1}^m\langle G_n^{PV},e_k\rangle .
$$ 
In the following, when no confusion is possible, we lighten the notation $G_n^{PV}$ into $G_n$, and $L_n^{V}$ into $L_n$, like in Theorems \ref{asymptoticsDTG} and \ref{asymptoticsDTL} below.

%%%%%%%%%%%%%%%%%%%%%%%%%
\section{Phase transition: \emph{small} and \emph{large} B-trees}
\label{sec-transition}

With regard to the asymptotics of their composition vector, P\'olya urns are subject to a well known phase transition. See the above references on P\'olya urns for a general treatment. Let us translate the phenomenon for our B-urns, looking at 
 the spectral properties of the replacement matrix. 

%%%%%%%%%%%%%
\subsection{Spectral decomposition of the ambient space}

In this section, we state notations relative to the spectral decomposition of the matrix~$R_m$.
These notations are used all along the paper.
The (unitary) characteristic polynomial of $R_m$ turns out to be
\begin{equation}
\label{polycar}
\chi _m(X)=
\prod _{k=m}^{2m-1}(X+k)-\frac{(2m)!}{m!}.
%=\frac{(2m)!}{m!}\left( \prod _{k=m}^{2m-1}\frac{X+k}{1+k}-1\right)
\end{equation}
Its complex roots are all simple, the one having the largest real part one being $1$.
Furthermore, two distinct eigenvalues that have the same real part are conjugated.
We denote by
\begin{equation}
\label{deflambda2}
\lambda _2=\sigma _2+i\tau _2
\end{equation}
the eigenvalue of $R_m$ having the second largest real part $\sigma _2$ and a positive imaginary part $\tau _2$.
We adopt also the following notations:
\begin{equation}
\label{notationsSpectrales}
\left\{
\begin{array}{l}
H_{m+1}(X)=\displaystyle
\sum _{k=1}^m\frac{1}{X+k}
\\ \\
v(\lambda )=\displaystyle
\frac{1}{(m+\lambda )H_{m+1}(m+\lambda -1)}\times
\\ \\
\hskip 20pt\displaystyle
\left(
1,\frac{m+1}{m+1+\lambda}, \frac{(m+1)(m+2)}{(m+1+\lambda)(m+2+\lambda)}, \dots ,\frac{(m+1)\dots(2m-1)}{(m+1+\lambda)\dots (2m-1+\lambda)}
\right)
\\ \\
\langle v(\lambda ),e_k\rangle = \displaystyle
\frac{1}{(m+\lambda )H_{m+1}(m+\lambda -1)} \prod_{j=1}^{k-1}\frac{m+j}{m+j+\lambda}
\\ \\
u(\lambda )\left( x_1,\dots ,x_m\right) =\displaystyle
\sum _{k=1}^m\left( \prod _{j=0}^{k-2}\frac{\lambda+m+j}{1+m+j}\right) x_k
\\ \\
\hskip 90pt  =x_1 + \displaystyle \frac{\lambda+m}{1+m}\ x_2 +\frac{(\lambda+m)(\lambda+m+1)}{(1+m)(2+m)}\ x_3 + \cdots
\\ \\
\hskip 100pt + \displaystyle\frac{(\lambda+m)\dots(\lambda+2m-2)}{(1+m)\dots(2m-1)}\ x_m.
\end{array}
\right.
\end{equation}
When $\lambda$ is an eigenvalue of $R_m$, the vector $v(\lambda )$ is an eigenvector of $\transp{R_m}$ associated with $\lambda$.
The linear form $u(\lambda )$ is an eigenform of $\transp{R_m}$ associated with $\lambda$, which means that for any (column) vector $V$, $u(\lambda )\left(\transp{R_m}V\right) =\lambda u(\lambda )(V)$.
Moreover, if $\lambda$ and $\mu$ are eigenvalues of $R_m$, then $u(\lambda)\left[ v(\mu )\right]=\delta _{\lambda ,\mu}$ (Kronecker).
In other words, $\left( v(\lambda )\right) _{\lambda\in\Sp\left( R_m\right)}$ and $\left( u(\lambda )\right) _{\lambda\in\Sp\left( R_m\right)}$ are dual basis of respectively eigenvectors and eigenforms of $\transp{R_m}$.

In the sequel, for more simplicity, we denote
\begin{equation}
\label{notations12}
\left\{
\begin{array}{l}
v_1=v(1)=\displaystyle
\frac 1{H_{m+1}(m)}\left( \frac 1{m+1}, \frac 1{m+2}, \dots , \frac 1{2m} \right)
\\ \\
u_1\left( x_1,\dots ,x_m\right) =u(1)\left( x_1,\dots ,x_m\right) =\displaystyle
\sum _{k=1}^mx_k
\\ \\
v_2=v(\lambda _2)
{\rm ~~and~~}
u_2=u(\lambda _2).
\end{array}
\right.
\end{equation}

The complex vector space $\g C^m$ admits the decomposition as direct sum of $\transp{R_m}$-stable lines
$$
\g C^m=\bigoplus _{\lambda\in\Sp\left( R_m\right)}\g Cv(\lambda )
$$
and the corresponding projection on any line $\g Cv(\lambda )$ is $u(\lambda )v(\lambda)$.
In the real field, we use the decomposition
\begin{equation}
\label{defRondV1}
\g R^m=\g Rv_1\oplus\rond V_1
\end{equation}
where $\rond V_1$ is the only subspace which is simultaneously $\transp{R_m}$-stable and supplementary to $\g Rv_1$.
It is generated by the vectors respectively constituted by the real parts and the imaginary parts of the coordinates of the complex vectors $v(\lambda )$, $\lambda\in\Sp\left( R_m\right)\setminus\{ 1\}$.
In the same vein, we denote by $\rond V_2$ the only $\transp{R_m}$-stable subspace of $\g R^m$ that satisfies
$$
\g R^m=\g Rv_1\oplus\g R\Re\left( v_2\right)\oplus\g R\Im\left( v_2\right)\oplus\rond V_2.
$$

%%%%%%%%%%%%%
\subsection{Phase transition for gaps and fringe nodes}

The phase transition on urns is expressed on the gap process $(G_n)_n$ in the following result.
Note the two very different convergence modes:
a weak one for small phases \emph{vs} a strong one with periodic phenomena for large phases.
See simulations in Section~\ref{sec-simulations} for an illustration.

\begin{Th}
\label{asymptoticsDTG}
Let $m\geq 2$.
Let $V\in\g N^m$ be a non zero vector and let $\left( G_n\right) _{n\geq 0}$ be the (discrete-time) gap process starting with the initial condition $G_0=PV$.
Then, with notations~\eqref{notationsSpectrales} and \eqref{notations12},

\vskip 5pt
(i) (Small phases)

if $m\leq 59$, as $n$ tends to infinity, $\displaystyle\frac{G_n-nv_1}{\sqrt n}$ converges in distribution to a centered Gaussian vector;

\vskip 5pt
(ii) (Large phases)

if $m\geq 60$, as $n$ tends to infinity,
\begin{equation}
\label{expansionDTG}
G_n=nv_1+2\Re \left( n^{\lambda _2}W^{DT}v_2\right) +o\left( n^{\sigma _2}\right),
\end{equation}
almost surely and in any ${\rm L}^p$, $p\geq 1$, where $W^{DT}$ is a complex-valued random variable with expectation $\displaystyle\frac{\Gamma\left( |PV|\right)}{\Gamma\left( |PV|+\lambda _2\right)}u_2(PV)$
(\hskip 2pt$\Gamma$ denotes Euler Gamma function).
\end{Th}

\pff

These results come from the general theory of balanced P\'olya urn processes.
See Janson~\cite{Jan} or~\cite{Pou08}. The numerical values of $\sigma_2$ leading to the phase transition are given in the Appendix.
$\hfill\fin$

\begin{Rem}
In Theorem~\ref{asymptoticsDTG}(ii), for any $p\in ]0,1[$, the asymptotics is also valid in the (not locally convex) complete metric space ${\rm L}^p$ defined by the usual quasi-norm.
This is true for Theorem~\ref{asymptoticsDTL}, Corollary~\ref{asymptoticsBtreeL}, Theorem~\ref{asymptoticsCTG}, Theorem~\ref{asymptoticsCTL} and Remark~\ref{rem-expansionBtreeCT} as well.
\end{Rem}

When $V$ is a complex vector, $\Re (V)$ denotes the vector made of real parts of $V$'s coordinates.
The random variable $W^{DT}$, which is more deeply studied below, appears as a martingale limit in the field of urn theory.
It can be also described as the almost sure limit of $G_n$ after normalisation and projection along the principal direction defined by $v_2$:
$$
W^{DT}=\lim _{n\to\infty}\frac{1}{n ^{\lambda _2}}u_2\left( G_n\right).
$$
Of course, the phase transition and the asymptotics can be straightforwardly translated on the fringe node process $\left( L_n\right) _n$ \emph{via} the diagonal matrix $P$ defined in \eqref{defP}.
The random variable $W^{DT}$ that appears for large phases in Theorem~\ref{asymptoticsDTL} has the same law as the one of the gap process in Theorem~\ref{asymptoticsDTG}.

\begin{Th}
\label{asymptoticsDTL}
Let $m\geq 2$.
Let $V\in\g N^m$ be a non zero vector and let $\left( L_n\right) _{n\geq 0}$ be the (discrete-time) fringe node process starting with the initial condition $L_0=V$.
Then, with notations~\eqref{notationsSpectrales}, \eqref{notations12} and~\eqref{defP},

\vskip 5pt
(i) (Small phases)

if $m\leq 59$, as $n$ tends to infinity, $\displaystyle\frac{L_n-nP^{-1}v_1}{\sqrt n}$ converges in distribution to a centered Gaussian vector;

\vskip 5pt
(ii) (Large phases)

if $m\geq 60$, as $n$ tends to infinity,
\begin{equation}
\label{expansionDTL}
L_n=nP^{-1}v_1+2\Re \left( n^{\lambda _2}W^{DT}P ^{-1}v_2\right) +o\left( n^{\sigma _2}\right),
\end{equation}
almost surely and in any ${\rm L}^p$, $p\geq 1$, where $W^{DT}$ is a complex-valued random variable which has the same distribution as in the variable named the same way in Theorem~\ref{asymptoticsDTG}.
\end{Th}

\vskip 15pt

Geometrically speaking, expansion \eqref{expansionDTG} (and expansion \eqref{expansionDTL} as well) can be understood as follows. Notice that an analogous explanation holds for expansions \eqref{expansionCTG} and \eqref{expansionCTL} in Theorem \ref{asymptoticsCTG} and Theorem \ref{asymptoticsCTL} respectively.

\vskip 5mm

\begin{figure}[!h]
\label{spirale}
\begin{center}
\begin{picture}(420,190)
\put(217,-118){\includegraphics[height=400pt,width=250pt]{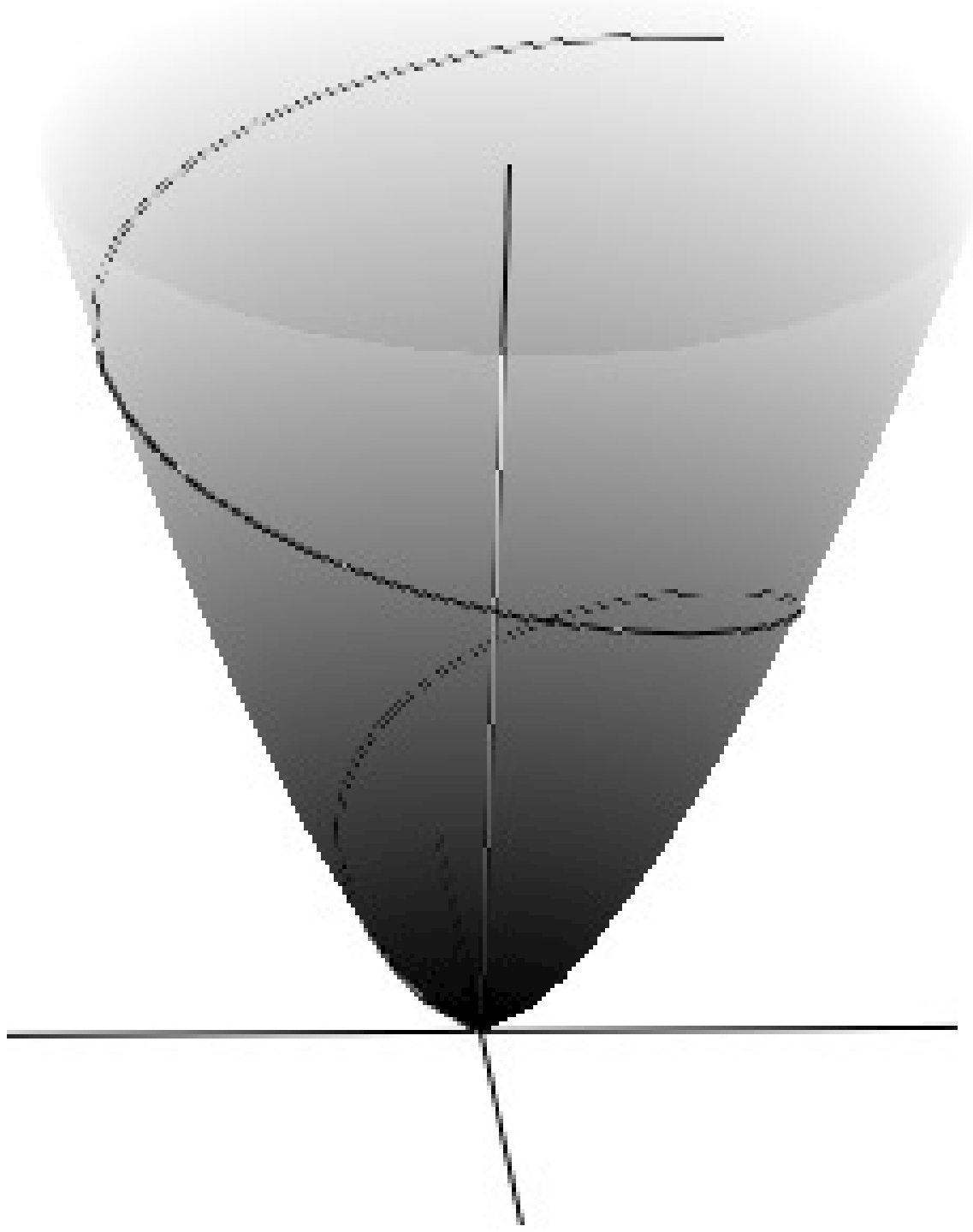}}
\put(350,-10){$\Re (v_2)$}
\put(380,25){$\Im (v_2)$}
\put(350,145){$v_1$}
\put(0,80){\begin{minipage}{250pt}
Let us denote by $\varphi$ any argument of the complex number $W^{DT}$.
T he trajectory of the random vector $G_n$, projected in the $3$-dimensional real vector space spanned by the vectors $(\Re (v_2),\Im (v_2),v_1)$ is almost surely
asymptotic to the (random) spiral
$$
\left\{
\begin{array}{l}
x_n=2|W|n^{\sigma _2}\cos (\tau _2 \log n+\varphi), \\
y_n=-2|W|n^{\sigma _2}\sin (\tau _2 \log n+\varphi), \\
z_n=n,
\end{array}
\right.
$$
drawn on the (random) revolution surface
$$
4|W|^2z^{2\sigma _2}=x^2+y^2,
$$
when $n$ tends to infinity.
\end{minipage}}
\end{picture}
\end{center}
%\caption{}
\end{figure}

\vskip 10pt
As is well known in the field of P\'olya urn processes, the phase transition is due to the number $\sigma _2$.
When $\sigma _2<1/2$, the P\'olya urn is \emph{small} and admits a weak Gaussian asymptotics.
On the contrary, when $\sigma _2>1/2$, the urn is \emph{large} and has a strong and oscillating ($\lambda _2$ is nonreal) asymptotic behaviour.
Considering the replacement matrix $R_m$, it turns out that $\sigma _2$ is an increasing function of $m$ and that:

$\bullet$
when $m=59$, $\lambda _2=(0.49534...)+(9.10305...)i$

while

$\bullet$
when $m=60$, $\lambda _2=(0.50378...)+(9.10270...)i$.

These numerical values have been computed by a Newton approximation algorithm, which can be found in the Appendix. The monotonicity of $\sigma _2$ as a function of $m$ (it increases to $1$ when $m$ tends to infinity) has been evocated by Hennequin \cite{Hen} in a figure. Qualitatively, let us emphasize the fact that for large values of $m$ (which is the actual use in computer science, since $m$ amounts to several hundreds), the fluctuation term with $W^{DT}$ is highly significant.

\vskip 10pt
We deduce from these theorems the asymptotic behaviour of the composition vector of the fringe nodes of different types in a B-tree, which is a particular case of Theorem~\ref{asymptoticsDTL} with the initial condition $V=(1,0,\dots ,0)$.
We stated the theorems above for arbitrary initial conditions because of the further study of the limit law $W^{DT}$ that requires these wider statements.

\begin{Cor}
\label{asymptoticsBtreeL}
Let $m\geq 2$.
Let $\rond L_n$ be the composition vector at time $n$ of the fringe nodes of different types in a B-tree with minimum degree $m$.
Then, as $n$ goes of to infinity, with notations~\eqref{defP} and~\eqref{notations12},

\vskip 5pt
(i)
when $m\leq 59$, $\displaystyle\frac{\rond L_n-nP^{-1}v_1}{\sqrt n}$ converges in distribution to a centered gaussian vector;

\vskip 5pt
(ii)
when $m\geq 60$, $\rond L_n=nP^{-1}v_1+2\Re \left( n^{\lambda _2}W_m^{\rm B-tree}P ^{-1}v_2\right) +o\left( n^{\sigma _2}\right)$, almost surely and in any ${\rm L}^p$, $p\geq 1$, where $W_m^{\rm B-tree}$ is a complex-valued random variable with expectation $\displaystyle\frac{m!}{\Gamma\left( m+\lambda _2\right)}$.
\end{Cor}

%%%%%%%%%%%
%\subsection{Simulations illustrating the phase transition}

%On B-trees, ie on $\rond L_n$ and on $W_m^{\rm B-tree}$.

%%%%%%%%%%%%%%%%%%%%%%%%%
\section{Embeddings into continuous time}
\label{sec-embedding}

Going further and obtaining significant properties of the random limit $W^{DT}$ is not so easy. As can be seen in this section, the classical method of embedding in continuous time turns out to be very fruitful: 
this idea of embedding discrete urn models in continuous-time branching processes goes back
at least to Athreya and Karlin  \cite{AK}.
A description is given in Athreya and Ney's book \cite[Section 9]{AN}.
The method has been recently revisited and developed by Janson \cite{Jan}, it is the core of recent results on P\'olya urns in \cite{ChaPouSah, ChaMaiPou}.

%%%%%%%%%%%%%%
\subsection{Definition of the continuous-time fringe node process}

Denote by $\left( L(t)\right) _{t\in\g R_{\geq 0}}$ the $\g N^m\setminus\{ 0\}$-valued continuous time Markov process having $\rond G$ as infinitesimal generator, where $\rond G$ is defined, for any function $f:\g N^m\setminus\{ 0\}\to \rond V$ ($\rond V$ is any real or complex vector space) and for any nonzero $X=\left( x_1,\dots ,x_m\right)\in\g N^m$, by
$$
\rond G(f)(X)=
\sum _{k=1}^m(m+k-1)x_k\Big[ f\left( X+w_k\right)- f\left( X\right)\Big]
$$
where the increment vectors $w_k$ have already be defined by~\eqref{increments}.

\vskip 10pt
This process is a multitype branching process, embedding of the Markov chain $\left( L_n\right) _n$ into continuous time, as classically done (see for example Bertoin \cite{Bertoin}).
One can think of it the following way.
At each (real) time $t\geq 0$, one gets particles of $m$ different types named $1,2,\dots ,m$.
Each particle is equipped with a clock that rings at random times.
The clock of any particle of type $k$ is exponentially distributed, with parameter $m+k-1$ and all the clocks are independent.
The dynamics of the process is the same as in discrete time:
for any $k\in\{ 1,\dots ,m-1\}$, when the clock of a particle of type $k$ rings, the particle disappears and is replaced by a particle of type $k+1$;
when the clock of a particle of type $m$ rings, the particle disappears and is replaced by two particles of type $1$.

\vskip 10pt
Having the same dynamics, the distributions of the processes $\left( L_n\right) _{n\in\g N}$ and $\left( L(t)\right) _{t\in\g R_{\geq 0}}$ are as usual related by the finite-time connection
\begin{equation}
\label{finiteTimeConnectionL}
(L_n)_{n\in\g N}\egalLoi \left(L\left(\tau _{(n)}\right)\right)_{n\in\g N}
\end{equation}
where $\tau _{(n)}$ denotes the $n$-th splitting time
(the $n$-th ringing time).
This relation allows us to transfer results on one process to the other. In particular, the results below strongly rely on the fact that $(e^{- (\transp{R_m}) t}L(t))_{t\in\g R_{\geq 0}}$ is a vector-valued martingale (see Janson \cite{Jan} or Athreya and Ney \cite{AN} for this).

\subsection{Definition of the continuous-time gap process}

Define the vector-valued continuous-time Markov process $\left( G(t)\right) _{t\in\g R_{\geq 0}}$ as being the embedding into continuous time of the discrete-time urn process $\left( G_n\right) _{n\in\g N}$.
It takes its values in the set of vectors of the form $PV$ where $V\in\g N^m\setminus\{ 0\}$.
With notations as above, its infinitesimal generator is given by
$$
\rond H(f)(X)=
\sum _{k=1}^mx_k\Big[ f\left( X+Pw_k\right)- f\left( X\right)\Big] ,
$$
the increment vectors $Pw_k$ being the rows of the urn replacement matrix $R_m$.

\vskip 10pt
One can think of this process the following way.
Take an urn that contains clocks of $m$ different colors named $1,\dots ,m$.
Each clock rings at a random time, exponentially distributed with parameter $1$ and all the clocks are independent.
As soon as a clock rings, the following replacement mechanism occurs:
if the ringing clock has color $k\in\{ 1,\dots ,m-1\}$, then it disappears together with $m+k-2$ other clocks of color $k$ and $m+k$ clocks of color $k+1$ arise in the urn;
if the ringing clock has color $m$, then it disappears together with $2m-2$ other clocks of color $k$ and $2m$ clocks of color $1$ arise in the urn.
The fact that the ringing times are exponentially distributed allows to think as if all clocks were restarted as soon as one of them rings.
Note that the fact that many clocks disappear at the same time prevents $\left( G(t)\right) _{t\in\g R_{\geq 0}}$ from being a multitype branching process.

\vskip 10pt
As in the preceding case, the processes $\left( G_n\right) _{n\in\g N}$ and $\left( G(t)\right) _{t\in\g R_{\geq 0}}$ have the same dynamics, so that
\begin{equation}
\label{finiteTimeConnectionG}
\left(G_n\right)_{n\in\g N}\egalLoi \left(G\left(\tau '_{(n)}\right)\right)_{n\in\g N},
\end{equation}
where $\tau '_{(n)}$ denotes the $n$-th ringing time.

\vskip 10pt
When $V\in\g N^m\setminus\{ 0\}$, denote by $\left( L(t)^V\right) _{t\geq 0}$ the fringe node process starting with $L(0)=V$ and by $\left( G(t)^V\right) _{t\geq 0}$ the gap process starting from $G(0)=V$.
Then, as in the discrete-time case in \eqref{lienLeafGapDiscret}, for any $V\in\g N^m\setminus\{ 0\}$,
\begin{equation}
\label{lienLeafGapContinu}
\left( G(t)^{PV}\right) _{t\geq 0}
\egalLoi
\left( P L(t)^V\right) _{t\geq 0},
\end{equation}
where $P$ is the diagonal matrix defined in \eqref{defP}.

\subsection{Asymptotics of both continuous-time processes}

The asymptotics of the continuous-time processes admit the same kind of phase transition as in discrete time.
We state this asymptotics for both continuous-time processes in Theorems~\ref{asymptoticsCTL} and~\ref{asymptoticsCTG}.
Since $G$ is the image of $L$ by $P$, any of these theorem implies the other one.
Nevertheless, as explained below, we prove both of them together using results on branching processes and results on P\'olya urns.

\begin{Th}
\label{asymptoticsCTG}
Let $m\geq 2$.
Let $V\in\g N^m$ be a non zero vector and let $\left( G(t)\right) _{t\in\g R_{\geq 0}}$ be the continous-time gap process that satisfies $G(0)=PV$.
Then, with notations~\eqref{notationsSpectrales},

\vskip 5pt
(i) (Small phases)

when $m\leq 59$, as $t$ tends to infinity, $e^{-t}G(t)$ converges almost surely and in any ${\rm L}^p$, $p\geq 1$, to $\xi v_1$ where $\xi$ is a positive random variable which is Gamma-distributed with parameter $|PV|$.
Furthermore, if one writes $G(t)=G_1(t)+G'_1(t)$ where the random vector $G_1(t)$ is proportional to $v_1$ and where $G'_1(t)$ is $\rond V_1$-valued (see\eqref{defRondV1}), then $e^{-t}G_1(t)$ converges almost surely and in any ${\rm L}^p$ to $\xi v_1$ while $e^{-t/2}G'_1(t)$ converges in distribution to $\sqrt\xi N$ where $N$ is a centered $\rond V_1$-valued gaussian vector independant of $\xi$.

\vskip 5pt
(ii) (Large phases)

when $m\geq 60$, as $t$ tends to infinity,
\begin{equation}
\label{expansionCTG}
G(t)=e^t\xi v_1 \left( 1+o(1)\right)
+2\Re\left( e^{\lambda _2t}W^{CT}v_2\right)\left( 1+o(1)\right)
+o\left( e^{\sigma _2t}\right)
\end{equation}
almost surely and in any ${\rm L}^p$, $p\geq 1$, where $W^{CT}$ is a complex-valued random variable with expectation $u_2(PV)$ and $\xi$ a positive random variable that is Gamma distributed with parameter $|PV|$.
The almost sure remainder $o\left( e^{\sigma t}\right)$ is a $\rond V_2$-valued random vector.
\end{Th}

\begin{Th}
\label{asymptoticsCTL}
Let $m\geq 2$.
Let $V\in\g N^m$ be a non zero vector and let $\left( L(t)\right) _{t\in\g R_{\geq 0}}$ be the continous-time fringe node process that satisfies $L(0)=V$.
Then, with notations~\eqref{notationsSpectrales} and~\eqref{defP},

\vskip 5pt
(i) (Small phases)

when $m\leq 59$, as $t$ tends to infinity, $e^{-t}L(t)$ converges almost surely and in any ${\rm L}^p$, $p\geq 1$, to $\xi P^{-1}v_1$ where $\xi$ is a positive random variable which is Gamma-distributed with parameter $|PV|$.
Furthermore, if one writes $L(t)=L_1(t)+L'_1(t)$ where the random vector $L_1(t)$ is proportional to $P^{-1}v_1$ and where $L'_1(t)$ is $P^{-1}\rond V_1$-valued (see\eqref{defRondV1}), then $e^{-t}L_1(t)$ converges almost surely and in any ${\rm L}^p$ to $\xi P^{-1}v_1$ while $e^{-t/2}L'_1(t)$ converges in distribution to $\sqrt\xi N'$ where $N'$ is a centered $P^{-1}\rond V_1$-valued gaussian vector independant of $\xi$.

\vskip 5pt
(ii) (Large phases)

when $m\geq 60$, as $t$ tends to infinity,
\begin{equation}
\label{expansionCTL}
L(t)=e^t\xi P^{-1}v_1 \left( 1+o(1)\right)
+2\Re\left( e^{\lambda _2t}W^{CT}P^{-1}v_2\right)\left( 1+o(1)\right)
+o\left( e^{\sigma _2t}\right)
\end{equation}
almost surely and in any ${\rm L}^p$, $p\geq 1$, where $W^{CT}$ is a complex-valued random variable with expectation $u_2(PV)$ and $\xi$ a positive random variable that is Gamma distributed with parameter $|PV|$.
The almost sure remainder $o\left( e^{\sigma t}\right)$ is a $P^{-1}\rond V_2$-valued random vector.
\end{Th}

Note that the random variables $\xi$ and $W^{CT}$ that appear in both theorems have been denoted the same way because their distributions are the same in both cases.
This comes immediately from~\eqref{lienLeafGapContinu}.

\vskip 5pt
{\sc Proof of Theorems~\ref{asymptoticsCTG} and~\ref{asymptoticsCTL}}.

Despite the fact that similar results can be found in Janson \cite{Jan} and Mailler \cite{Mailler}, the particular case of our processes is not properly contained in their statement. The proofs are essentially made the same way as in both papers~\cite{ChaLiuPouContinu} and~\cite{ChaLiuPouDiscret}.
We give hereunder the general scheme of the argumentation.
The first tool comes from the fact that the normalised projection $\left( e^{-t}u_1(G(t))\right) _{t\geq 0}$ is always a convergent positive martingale.
The random variable $\xi$ is its limit.

\vskip 3pt
(i) Small phases.
The process $\left( L(t)\right) _{t\in\g R_{\geq 0}}$ is a multitype branching process so that {\it (i)} in Theorem~\ref{asymptoticsCTL} is covered by~\cite{AN} and~\cite{Jan}.
Relation~\eqref{lienLeafGapContinu} thus implies {\it (i)} in Theorem~\ref{asymptoticsCTG}.

\vskip 3pt
(ii) Large phases.
As for the first projection, $\left( e^{-\lambda _2t}u_2(G(t))\right) _{t\geq 0}$ is a martingale, which is convergent if, and only if $\sigma _2>1/2$, \emph{i.e.} when $m\geq 60$.
The complex-valued random variable $W^{CT}$ is its limit.
The oscillating term $\Re\left( e^{\lambda _2t}W^{CT}P^{-1}v_2\right)$ in Theorem~\ref{asymptoticsCTL} is a consequence of~\cite{AN} and~\cite{Jan}'s results.
In order to establish the almost sure remainders $o\left( e^{\sigma _2t}\right)$, we use results on discrete-time P\'olya urns shown in~\cite{Pou08}.
The work is done on the gap process $(G(t))_t$ viewed as an embedded urn into continuous time.
For any $t\geq 0$, decompose $G(t)$ as the sum $G(t)=G_1(t)+G_2(t)+G_\ell (t)+G_s(t)$ of its respective following projections on the described supplementary subspaces:

$\bullet$ $G_1(t)$ is the projection on $\g Rv_1$ as before;

$\bullet$ $G_2(t)$ is the projection on the real plane generated by the real part and the imaginary part of $v_2$;

$\bullet$ $G_\ell (t)$ is the projection on the subspace of $\g R^m$ generated by the real and imaginary parts of the eigenvectors $v(\lambda )$ for all eigenvalues $\lambda$ different from $1$ and $\lambda _2$ such that $\Re\left(\lambda\right) >1/2$
(\emph{large} projections);

$\bullet$  finally, $G_s(t)$ is the projection on the subspace of $\g R^m$ generated by the real and imaginary parts of the eigenvectors $v(\lambda )$ for all eigenvalues $\lambda$  such that $\Re\left(\lambda\right)\leq 1/2$
(\emph{small} projections).

As seen before, $e^{-t}G_1(t)$ converges to $\xi v_1$ almost surely and in ${\rm L}^p$, $p\geq 1$, by martingale techniques;
this gives rise to the first term $e^t\xi v_1$ in the asymptotics of $G(t)$.
Since $G_2(t)=2\Re\left[ u_2\left( G_2(t)\right)v_2\right]$ and because of the convergence in ${\rm L}^p$, $p\geq 1$, of the complex martingale $\left( e^{-\lambda _2t}u_2(G(t))\right) _{t\geq 0}$ mentioned above, one gets the second term $\Re\left( e^{\lambda _2t}W^{CT}v_2\right)$ of $G(t)$'s asymptotics.
The remainder $o\left( e^{\sigma _2t}\right)$ is obtained from $G_\ell$ and $G_s$ asymptotics.
As for $G_2$, if $\lambda$ is an eigenvalue of $R_m$ such that $\Re (\lambda )>1/2$, by martingale arguments, the complex projection of $G(t)$ on any eigenline $\g Cv(\lambda )$ is equivalent to $e^{\lambda t}W_\lambda$ almost surely and in any ${\rm L}^p$ where $W_\lambda$ is a complex-valued random variable.
In particular, the whole projection $G_\ell (t)$ is $o\left( e^{\sigma _2t}\right)$, almost surely and in any ${\rm L}^p$.

To make the proof complete, it remains to show that $G_s(t)$ is $o\left( e^{\sigma _2t}\right)$ as well.
To prove this fact, we use the technique detailed in~\cite{ChaLiuPouContinu} (Theorem 4.1 and Lemma 4.2).
It consists in considering the same projection for the discrete-time urn process $(G_n)_n$, in using the moment bounds proven in~\cite{Pou08} for small projections of discrete-time P\'olya urns and in coming back to continuous time by Relation~\eqref{finiteTimeConnectionG}.
By this means, one shows after some probabilistic arguments that for any $\eta >0$, 
the whole projection $G_s$ satisfies that $e^{-\left(\eta +\frac12\right)t}G_s(t)$ is bounded, almost surely and in ${\rm L}^p$, $p\geq 1$, implying the expected result on $G(t)$.
The corresponding asymptotics of $L(t)$ is obtained by taking the image of $G(t)$ by $P^{-1}$.
$\hfill\fin$

\medskip

\begin{Rem}
\label{rem-expansionBtreeCT}
For $m\geq 60$, we deduce from these theorems the asymptotic behaviour of the continuous-time fringe node process, denoted by $(\rond L(t))_t$, starting from the B-tree initial condition $V=(1, 0, \dots , 0)$:
\begin{equation}
\label{expansionBtreeCTL}
\rond L(t)=e^t\xi P^{-1}v_1 \left( 1+o(1)\right)
+2\Re\left( e^{\lambda _2t}\rond W^{CT}P^{-1}v_2\right)\left( 1+o(1)\right)
+o\left( e^{\sigma _2t}\right)
\end{equation}
almost surely and in any ${\rm L}^p$, $p\geq 1$, where $\rond W^{CT}$ is a complex-valued random variable with expectation $m$ and $\xi$ a positive random variable that is Gamma distributed with parameter $m$.
The almost sure remainder $o\left( e^{\sigma t}\right)$ is a $P^{-1}\rond V_2$-valued random vector.
\end{Rem}

\medskip

For large phases, the finite time connections~\eqref{finiteTimeConnectionG} or~\eqref{finiteTimeConnectionL} lead to a relation between the random variables $W$ in discrete and continuous times.
This relation, commonly named \emph{martingale connection} will be stated and used  below in the article.
We indicate hereafter how one can get it.
Take for instance Relation~\eqref{finiteTimeConnectionG} concerning the gap processes $(G_n)_n$ and $(G(t))_t$ starting with the same initial condition $G_0=G(0)=PV$.
Using Theorems~\ref{asymptoticsDTG} and~\ref{asymptoticsCTG}, since $\tau _{(n)}$ tends almost surely to $+\infty$ as $n$ goes of to infinity, one gets successively
$\xi =\lim _{t\to\infty}e^{-t}u_1\left( G(t)\right)
=\lim _{n\to\infty}e^{-\tau _{(n)}}u_1(G_n)
=\lim _{n\to\infty}ne^{-\tau _{(n)}}$
on one hand.
On the other hand,
$W^{CT}=\lim _{t\to\infty}e^{-\lambda _2t}u_2\left( G(t)\right)
=\lim _{n\to\infty}e^{-\lambda _2\tau _{(n)}}u_2(G_n)
=\lim _{n\to\infty}\left[ ne^{-\tau _{(n)}}\right]^{\lambda _2}\left[ n^{-\lambda _2}u_2(G_n)\right]$.
This entails the martingale connection
\begin{equation}
\label{martingaleConnection}
W^{CT}\egalLoi\xi ^{\lambda _2}W^{DT}.
\end{equation}
We just recall here that the random variable $\xi$ is Gamma-distributed with expectation $|PV|$.

%%%%%%%%%%%%%%%%%%%%%%%%%
\section{Limit law of large B-trees}
\label{sec-limitlaw}

In this section appear the benefits of the embedding in continuous time. Indeed, the branching property applied to the fringe node process $(L(t))_t$, together with the asymptotics proved in Theorem~\ref{asymptoticsCTL}, allow us to see the limit $W^{CT}$ as a solution of a  distributional equation. This is detailed in Section~\ref{sec-dislocation}. It is the starting point to deduce several properties of $W^{CT}$: its distribution is the unique solution of such an equation in a convenient space of probability distributions (Theorem \ref{th-contraction} in Section~\ref{sec-smoothing}); it admits exponential moments in a neighborhood of $0$ (Theorem \ref{Expmoments} in Section~\ref{sec-cascade}); up to a change of function, its Laplace transform is a solution of the quite simple (but unsolvable!) differential equation $y^{(m)}=y^2$ (Theorem \ref{equaDiffLaplaceCT} in Section~\ref{sec-laplace}); it admits a density relatively to Lebesgue measure on $\g C$ and its support is the whole complex plane (Theorem \ref{densiteSupportWCT} in Section~\ref{sec-density}). Thanks to connection \eqref{martingaleConnection} between $W^{CT}$ and $W^{DT}$,  corresponding results are true for $W^{DT}$ and consequently for $W_m^{\rm B-tree}$.

%%%%%%%%%%%%
\subsection{Dislocation equations in continuous time}
\label{sec-dislocation}

In this section, using the branching property of the continuous-time process $\left( L(t)\right)_t$, we show that the complex-valued random variable $W^{CT}$ is solution of a very simple distributional equation.
%As a consequence, we show that $W^{CT}$ admits a density and that it is supported by the whole complex plane.

In order to simplify the notations, for any $k\in\{ 1,\dots ,m\}$, denote by $W_k$ the limit random variable $W^{CT}$ (or its distribution) of the continuous-time fringe node process $\left( L(t)^{e_k}\right)_t$ that starts with one particle of type $k$, which means that its initial composition $L(0)$ is the $k$-th vector $e_k$ of $\g R^m$ canonical basis.
Denote also by $\tau _k$ the \emph{first} splitting time of the process $\left( L(t)^{e_k}\right)_t$;
its is exponentially distributed, with parameter $m+k-1$.

Because of the branching property of the process $\left( L(t)\right)_t$, for any time $t\geq \tau _1$, the processes $\left( L(t)^{e_1}\right)_{t\geq 0}$ and $\left( L(t)^{e_2}\right)_{t\geq 0}$ are related by the distributional equation 
$$L(t)^{e_1}\egalLoi L(t-\tau _1)^{e_2}.
$$
In the asymptotic form given by Theorem~\ref{asymptoticsCTL}, consider the second order term on both sides of the equality, which consists in projecting, normalizing and letting $t$ tend to infinity.
This leads to the distributional equality 
$$W_1=e^{-\lambda _2\tau _1}W_2,
$$ 
the random variables $W_2$ and $\tau _1$ being independent.
Doing the same for all values of $k\in\{1,\dots ,m\}$ leads to the distributional system:
\begin{equation}
\label{systLoiWCT}
\left\{
\begin{array}{ccl}
W_1&\egalLoi&e^{-\lambda _2\tau _1}W_2
\\
W_2&\egalLoi&e^{-\lambda _2\tau _2}W_3
\\
&\vdots&
\\
W_{m-1}&\egalLoi&e^{-\lambda _2\tau _{m-1}}W_m
\\
%W_m&\egalLoi&e^{-\lambda _2\tau _m}\left( W_1+W'_1\right)
W_m&\egalLoi&e^{-\lambda _2\tau _m}\left( W_1^{(1)}+W_1^{(2)}\right)
\end{array}
\right.
\end{equation}
where

$\bullet$
for any $k\in\{ 1,\dots ,m-1\}$, the random variables $\tau _k$ and $W_{k+1}$ of the $k$-th equation's right-hand sides are independent;

$\bullet$
in the right-hand side of the last equation, the random variables $W_1^{(1)}$ and $W_1^{(2)}$ are independent copies of $W_1$, both being independent of $\tau _m$ as well.

We recall that for any $k\in\{ 1,\dots ,m\}$, the random variable $\tau _k$ is exponentially distributed, with parameter $m+k-1$ (see Section~\ref{sec-embedding}).
In particular, $W_1$ is a solution of the following distributional equation, sometimes called fixed point equation or smoothing equation in some branching processes contexts (see Liu \cite{Liu98} or Biggins and Kyprianou \cite{BigKyp05}). 
\begin{equation}
\label{eqLoiWCT}
W_1\egalLoi B^{\lambda _2}\left( W_1^{(1)}+W_1^{(2)}\right)
\end{equation}
where

$\bullet$
the random variables $W_1^{(1)}$ and $W_1^{(2)}$ are independent copies of $W_1$;

$\bullet$
$B$ is a random variable, independent of $W_1^{(1)}$ and $W_1^{(2)}$, Beta distributed with parameters $(m,m)$ which means that it admits $t^{m-1}(1-t)^{m-1}\1 _{[0,1]}(t)$ as a density
($\1 _A$ denotes the indicatrix function of the set $A$).

The distribution of $B$ is computed the following way.
By immediate computation from System~\eqref{systLoiWCT}, one sees that $B=e^{-\left( \tau _1+\tau _2+\dots +\tau _m\right)}$, the variables $\tau _k$ being mutually independent.
To recognize the Beta$(m,m)$ law, one can make a direct computation of its density or compute its moments
(a Beta distribution is characterized by its moments because its support is compact).

%%%%%%%%%%%%
\subsection{Smoothing equation in discrete time}

In a general setting of $m$-color P\'olya urns, including the case of negative entries on the diagonal of the replacement matrix, Mailler \cite{Mailler} proves that $W^{DT}$ is a solution of a distributional equation which turns to be in our case
\begin{equation}
\label{eqLoiWDT}
W \ \egalenloi \  B_1^{ \lambda_2  }W^{(1)} + B_2^{ \lambda_2  }W^{(2)} ,
\end{equation}
where 

$\bullet$
the random variables $W^{(1)}$ and $W^{(2)}$ are independent copies of $W$;

$\bullet$
$(B_1,B_2)$ is a random vector, independent of $W^{(1)}$ and $W^{(2)}$, Dirichlet distributed with parameters $(m,m)$, which means that $B_1+B_2 = 1$ and that $B_1$ and $B_2$ are Beta distributed with parameters $(m,m)$.

\begin{Rem}
One proof of this result in \cite{Mailler} uses the tree structure of the urn. Nevertheless, we do not actually understand what kind of ``divide-and-conquer'' type argument, applied to B-trees, could lead to this equation. Indeed, in other cousin models, like $m$-ary search trees (see Fill and Kapur \cite{FillKapur}), such a backward decomposition leads to a finite time decomposition equation and passing to the limit, it gives the distributional equation. 
\end{Rem}

%%%%%%%%%%%%
\subsection{Contraction methods}
\label{sec-smoothing}

The question of existence and unicity of solutions of equations like \eqref{eqLoiWCT} or \eqref{eqLoiWDT} is classically solved using the Banach fixed point theorem. One point of view, frequent in analysis of algorithms, consists in starting from a decomposition property of the algorithm at finite time, deduce a distributional equation on a cost variable, and pass to the limit to get a smoothing equation on the limit random variable. See Knape and Neininger \cite{KnaNei} for P\'olya urns, and also the general paper by Neininger and R\"uschendorf \cite{NeiRusaap} or their survey \cite{NeiRusSurvey} for many examples of this so-called contraction method. Another point of view (in this article) consists in taking advantage of the dynamics of the algorithm and exhibiting a martingale limit, solution of a smoothing equation. Thus, the existence is automatically achieved. In both points of view, to get the unicity, the contraction property has to be established, in a convenient space of probability distributions, classically equipped with a Wasserstein distance to get a  complete metric space of measures.

We do not prove here the theorem below, since it is done in a general frame by Mailler \cite{Mailler}. The same kind of results can be found in
 Janson \cite[proof of Th 3.9 (iii)]{Jan} and in Knape and Neininger \cite{KnaNei} even if the only case $a_{i,i}\geq -1$ is considered there. See also \cite{ChaLiuPouContinu, ChaLiuPouDiscret}.

\begin{Th}
\label{th-contraction}
When $A$ is a complex number, let $\rond M_2\left( A\right)$ be the space of probability distributions on $\g C$ that have $A$ as expectation and a finite second moment, endowed with a complete metric space structure by the Wasserstein distance. Let $\lambda\in\C$ be any root of the characteristic polynomial~(\ref{polycar}) such that $\Re(\lambda) >\frac 12$. Then,
\begin{itemize}
\item[(i)]
Each of the two equations 
$$W \ \egalenloi \  B_1^{ \lambda }W^{(1)} + B_2^{ \lambda  }W^{(2)} 
$$ 
where $W^{(1)}$ and $W^{(2)}$ are independent copies of $W$, and where
$(B_1,B_2)$ is a random vector, independent of $W^{(1)}$ and $W^{(2)}$, Dirichlet distributed with parameters $(m,m)$, 
and 
$$W\egalLoi B^{\lambda }\left( W^{(1)}+W^{(2)}\right)
$$
where $W^{(1)}$ and $W^{(2)}$ are independent copies of $W$, and where
$B$ is independent of $W^{(1)}$ and $W^{(2)}$, Beta distributed with parameters $(m,m)$,

have a unique solution in $\rond M_2\left( A\right)$.
\item[(ii)]
For $m\geq 60$, the variable $W_m^{\rm B-tree}$, defined in Corollary \ref{asymptoticsBtreeL}, is the unique solution of \eqref{eqLoiWDT} having $\displaystyle\frac{m!}{\Gamma\left( m+\lambda _2\right)}$ as expectation and a finite second moment.
\item[(iii)]
For $m\geq 60$, the variable $\rond W^{CT}$, defined in \eqref{expansionBtreeCTL}, is the unique solution of \eqref{eqLoiWCT} having $m$ as expectation and a finite second moment.
\end{itemize}
\end{Th}

%%%%%%%%%%%%
\subsection{Cascades and exponential moments}
\label{sec-cascade}

Let $\lambda\in\C$ be any root of the characteristic polynomial~(\ref{polycar}) such that $\Re(\lambda) >\frac 12$ and let $B$ be a Beta distribution with parameters $(m,m)$. 
A simple computation leads to $2\g E\left(B^{\lambda}\right)=1$. This is coherent with equation 
\begin{equation}
\label{eqLoiWCTgeneral}
W\egalLoi B^{\lambda }\left( W^{(1)}+W^{(2)}\right).
\end{equation}
Moreover, for any positive real $s$,
$$
2\g E\left( B^s\right) = \frac{(2m)\dots (m+1)}{(2m-1+s) \dots (m+s)}<1\Longleftrightarrow s>1.
$$
Consequently, when $2\Re (\lambda) >1$, one has 
\begin{equation}
\label{coeffcontraction}
2\g E\left( |B^{\lambda}|^2\right) <1.
\end{equation}
Theorem \ref{Expmoments} below states that any solution $W$ of Equation \eqref{eqLoiWCTgeneral} admits exponential moments in a neighbourhood of $0$, so that the moment exponential generating series of $W$ defines an analytic function in a neighbourhood of the origin. Another consequence is that the law of $W$ is determined by its moments.

The proof relies on a Mandelbrot's cascade here defined in a complex setting (see Barral et al. \cite{Ba10} for complex Mandelbrot's cascades).

To lighten the notations, denote for a while $A:=B^{\lambda}$ and let $A_u, u\in U$ be independent copies of $A$, indexed by all finite sequences of $0$ and $1$:
$$
u = u_1\dots u_n \in U:= \bigcup_{n\geq 1} \{ 0, 1\}^n.
$$
Let $Y_0=m$, $Y_1 = 2mA$ and for $n\geq 2$,
 \begin{equation*}
 Y_n =  \sum_{u_1\dots u_{n-1} \in \{0,1\}^{n-1}}  2mAA_{u_1} A_{u_1u_2} \dots A_{u_1 \dots  u_{n-1}}.
 \end{equation*}
By the branching property, and using $2\g EA=1$, it is easy to see that $(Y_n)_n$ is a martingale with expectation  $m$. This martingale has been studied by many authors in the real-valued random variable case, especially in the context of Mandelbrot's cascades,
see for example Liu $\cite{Liu01}$  and the references therein.
It can be easily seen that
\begin{equation}
\label{EqYn}
    Y_{n+1} = B^{\lambda} \left( Y_n^{(1)} + Y_n^{(2)} \right)
\end{equation}
where $Y_n^{(1)}$ and $Y_n^{(2)}$ are independent of each other and independent of $B^{\lambda}$ and  each has the same distribution as $Y_n$. Therefore for $n\geq 1$, $Y_n$ is square-integrable and
$$ \Var Y_{n+1} =  2\g E|B^{\lambda}|^2 \Var Y_n + (4\g E |B^{\lambda}|^2  - 1)  $$
where $\Var X = \g E\left( |X-\g EX|^2\right)$  denotes the variance of $X$. Since $2\g E|B^{\lambda}|^2 < 1$, the martingale  $(Y_n)_n$ is bounded  in $L^2$, so that the following result holds.
\begin{Lem}
\label{YnL2}
Let $\lambda\in\C$ be any root of the characteristic polynomial~(\ref{polycar}) such that $\Re(\lambda) >\frac 12$ and let $B$ be a Beta distribution with parameters $(m,m)$.  When $n\rightarrow + \infty$,
 \begin{equation*}
    Y_n \rightarrow  Y_\infty \mbox{ a.s. and in } L^2,
  \end{equation*}
where $Y_\infty$ is a (complex-valued) random variable with variance
 \begin{equation*}
  \Var (Y_\infty) = \frac {4\g E |B^{\lambda}|^2  - 1}{1  - 2\g E|B^{\lambda}|^2 }.
 \end{equation*}
\end{Lem}
Notice that, passing to the limit in~(\ref{EqYn}) gives a new proof of the existence of a solution $W$ of  Eq. (\ref{eqLoiWCTgeneral}) with a given expectation and finite second moment whenever $\Re (\lambda) >1/2$. From Section \ref{sec-smoothing}, we have the uniqueness of solution of this equation so that Theorem~\ref{Expmoments} below will be proved as soon as it holds for $Y_\infty$.

\begin{Lem}
\label{expmomentsY}
There exist some constants $C>0$ and $\varepsilon >0$ such that for
all $t \in \C$ with $|t| \leq \varepsilon $, we have
\begin{equation}
 \label{bound-phi-inf}
  \g E  e^{\langle t, Y_\infty \rangle  } \leq e^{m\Re (t)  + C |t|^2 }.
 \end{equation}
\end{Lem}
\pff
By Fatou lemma, it is sufficient to prove the existence of $C>0$ and $\varepsilon >0$ such that for all $t \in \C$ with $|t| \leq \varepsilon $, and for every integer $n$,
\begin{equation}
 \label{bound-phi-n}
  \g E  e^{\langle t, Y_n \rangle  } \leq e^{m\Re (t)  + C |t|^2 }.
 \end{equation}
 Denote $\varphi_n(t):=\g E  e^{\langle t, Y_n \rangle  }$ and notice that $\varphi_{n+1}(t) = \g E\left( \varphi_n^2(tB^{\overline{\lambda}})\right)$ thanks to Equation~\eqref{EqYn}, allowing to prove \eqref{bound-phi-n} by recursion on $n\geq 0$. For $n=0$, 
 $$\varphi_0(t):=\g E  e^{\langle t, Y_0 \rangle  } = e^{m\Re(t)}
$$
and by the recursion assumption,
$$
\varphi_{n+1}(t)\leq \g E\left( e^{2C|t|^2 |B^{\lambda}|^2 + 2m\Re\left( tB^{\overline{\lambda}}\right)}\right)= e^{m\Re (t)  + C |t|^2 }f(t_1,t_2)
$$
where for any $t \in \C$, written $t = t_1+it_2$ with $t_1,t_2\in\g R$, 
$$f(t_1,t_2) = \g E\left( e^{C|t|^2 (2|B^{\lambda}|^2-1) + 2m\Re\left( tB^{\overline{\lambda}}\right) -m\Re(t)}\right),
$$
so that it is sufficient to prove that $(0,0)$ is a local maximum of $f$. Writing $\lambda = \sigma + i\tau$, with $\sigma, \tau\in\g R$,
$$
f(t_1,t_2) = \g E\exp\left[ C(t_1^2 + t_2^2)(2B^{2\sigma}-1) + 2m B^{\sigma}(t_1\cos(\tau) + t_2\sin(\tau)) -mt_1\right].
$$
Remembering $2\g E\left(B^{\lambda}\right)=1$, which means $2\g E\left( B^{\sigma}\cos(\tau)\right)=1$ and $\g E\left( B^{\sigma}\sin(\tau)\right)=0$, we get that the first derivatives vanish at $(0,0)$ which is a critical point. Moreover, the calculation of the second partial derivatives gives
\begin{eqnarray*}
 \frac{\partial^2 f}{\partial t_1} (0,0) &=& \g E \left[ \left(2mB^{\sigma}\cos(\tau)-m\right)^2+ 2C \left(2B^{2\sigma}-1 \right)\right],  \\
 \frac{\partial^2 f}{\partial t_2} (0,0)&= &\g E \left[ \left(2mB^{\sigma}\sin(\tau)\right)^2+ 2C \left(2B^{2\sigma}-1 \right)\right] ,\\
 \frac{\partial^2 f}{\partial t_1 \partial t_2} (0,0) &=& \g E \left(2mB^{\sigma}\cos(\tau)-m\right)\left(2mB^{\sigma}\sin(\tau)\right).
\end{eqnarray*}
By \eqref{coeffcontraction}, $\g E\left(2B^{2\sigma}-1 \right) <0$, 
so that  the Hessian matrix  at $(0,0)$ is definite negative for $C >0$ large enough which implies that $(0,0)$ is a local maximum of $f$. 
\hfill$\fin$

\medskip

The following theorem is a direct consequence of Lemma~\ref{expmomentsY}, like in~\cite{ChaLiuPouContinu}.

\begin{Th}
\label{Expmoments}
Let $\lambda \in \g C$ be a root of the characteristic polynomial~(\ref{polycar}) with $\Re (\lambda) >1/2$ and let $W$ be a solution of Eq. (\ref{eqLoiWCTgeneral}).
There exist some constants $C>0$ and $\varepsilon >0$
such that for all $t\in \C$ with $ |t| \leq \varepsilon$,
 \begin{equation}
 \label{bdExpmoments}
  \g Ee^{ \langle t,W\rangle   } \leq      e^{ m\Re (t)  + C |t|^2 }  \; \mbox{ and } \;
  \g Ee^{ |tW | } \leq   4   e^{ m|t|  +  2C |t|^2 } .
 \end{equation}
 \end{Th}

%%%%%%%%%%%%
\subsection{Laplace transform}
\label{sec-laplace}

Theorem \ref{Expmoments} above concerning $W^{CT}$ and Theorem 7 (ii) in Mailler \cite{Mailler} (which establishes that the Laplace series of $W^{DT}$ has an infinite radius of convergence) answer the question of the convergence of the Laplace series of $W^{DT}$ and $W^{CT}$. Nevertheless, a natural investigation consists in searching more information about these Laplace transforms coming from the smoothing equations.

\medskip

Indeed, the dislocation equations \eqref{systLoiWCT} lead to a system of differential equations on the Laplace transforms
\begin{equation}
\label{defLaplaceWCT}
\forall k = 1, 2, \dots , m, \hskip 1cm \varphi_k(z):= \g E\left( e^{<z,W_k>}\right),
\end{equation}
where we recall that $W_k$ is the limit random variable $W^{CT}$ of the continuous-time fringe node process $\left( L(t)^{e_k}\right)_t$ that starts with one particle of type $k$.
Using the independence between the splitting times $\tau_k$ (which are exponentially distributed, with parameter $m+k-1$) and the $W_k$, for $k = 1, 2, \dots , m-1$,
$$
\varphi_k(z) = \int_0^{+\infty} \varphi_{k+1}\left(ze^{-\overline{\lambda_2} t}\right)(m+k-1)e^{-(m+k-1)t}dt,
$$
and after a change of variable, and derivation, for $k = 1, 2, \dots , m-1$,
$$
\frac{m+k-1}{\overline{\lambda_2}}\varphi_k(z) +z\varphi_k'(z) = \varphi_{k+1}(z),
$$
and for $k=m$,
$$
\frac{2m-1}{\overline{\lambda_2}}\varphi_m(z) +z\varphi_m'(z) = \varphi_{1}^2(z).
$$
Thanks to a convenient change of function, a simple calculation gives the following theorem about the Laplace transforms of the $W_k$, for $k = 1, 2, \dots , m$. \begin{Th}
\label{equaDiffLaplaceCT}
For $k = 1, 2, \dots , m$, let $\varphi_k$ be the Laplace transform of $W_k$ defined in \eqref{defLaplaceWCT}, and let 
$$
\psi_k(z):= \left( -\overline{\lambda_2} \right)^{m+k-1} \ \frac{\varphi_k\left( z^{ -\overline{\lambda_2} }\right)}{z^{m+k-1}},
$$
(for any determination of the logarithm).
Then the functions $\psi_k$ satisfy the simple differential system
\begin{equation*}
\left\{
\begin{array}{l}
\displaystyle
\psi'_k=\psi_{k+1}, \hskip 1cm \forall k\in\{1,\dots ,m-1\},
\\
\displaystyle
\psi'_{m}=\psi_1^2.
\end{array}
\right.
\end{equation*}
In particular, $\psi_1$ is a solution of the differential equation 
$$y^{(m)}=y^2.
$$
\end{Th}

%%%%%%%%%%%
\subsection{Density and support}
\label{sec-density}

Liu's method has been developped in \cite{Liu99} and \cite{Liu01} for positive real-valued random variables solution of smoothing equations of the same type as \eqref{eqLoiWCT} or \eqref{eqLoiWDT}. Adapting this method to $\g C$-valued random variables, Mailler \cite{Mailler} gets the support and the existence of a density for the limit law of a $d$-color P\'olya urn. The theorem below is a particular case.

\begin{Th}
\label{densiteSupportWCT}
Let $m\geq 2$.
Let $V\in\g N^m$ be a non zero vector and let $W^{DT}$ and $W^{CT}$ be the $W$-distributions of the respective discrete-time and continuous-time fringe node processes having $V$ as initial composition
(see (ii) in Theorem~\ref{asymptoticsDTL} and Theorem~\ref{asymptoticsCTL}).
Then

\vskip 3pt
(i)
the supports of $W^{DT}$ and $W^{CT}$ are both the whole complex plane~$\g C$;

\vskip 3pt
(ii)
$W^{DT}$ and $W^{CT}$ are absolutely continuous relatively to Lebesgue's measure on~$\g C$;

\vskip 3pt
(iii)
as $|t| \rightarrow \infty$, $ \g Ee^{i\langle t,W^{DT}\rangle } = O(|t|^{-a})$ for each $a\in\left] 0, \displaystyle\frac m{\Re (\lambda_2)}\right[$, and the same is true for the Fourier transform of $W^{CT}$ as well.

\end{Th}

%%%%%%%%%%%
\subsection{Perspectives}
\label{sec-perspectives}

Some open questions remain about $W^{DT}$ and $W^{CT}$, let us say $W$: 

%- what is the true order of magnitude of $o(n^{\sigma_2})$ in \eqref{expansionDTG} or \eqref{expansionDTL}?

- can the $W$ distribution be expressed by means of usual distributions? Same question for $|W|$ and $\Arg(W)$?

- how heavy are the tails of $W$?

- what is the order of magnitude of $W$'s $p$-th moment as $p$ tends to $+\infty$?

%%%%%%%%%%%%%%%%
\section{Simulations}
\label{sec-simulations}

Let us here summarize and illustrate the asymptotic results concerning the fringe of  a random B-trees.
We only show the behaviour of the gap process $(G_n)_n$ and illustrate Theorem~\ref{asymptoticsDTG}; nevertheless, the simulations would be analogous for the fringe node process $(\rond L_n)_n$ to illustrate Corollary~\ref{asymptoticsBtreeL}.

In all the simulations below, sequences of $10^7$ random keys have been inserted in a B-tree for different value of the parameter $m$.
Notice that in ``real life'' computer science implementations, $m$ is most of the time taken around $100$ or more.

%%%%%%%%%%
\subsection{Simulations of $G_n$}
\label{sec-simulationsGn}

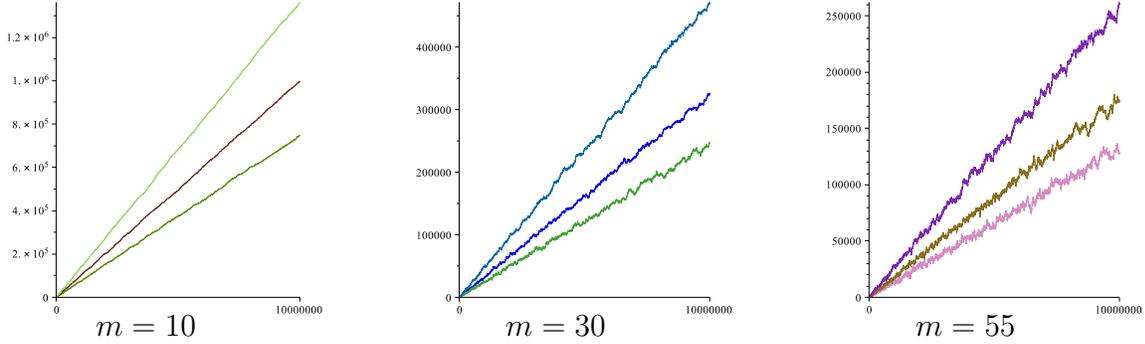
\begin{figure}[h]
\begin{center}
\begin{picture}(600,150)
\put(-10,0){\includegraphics[height=45 truemm]{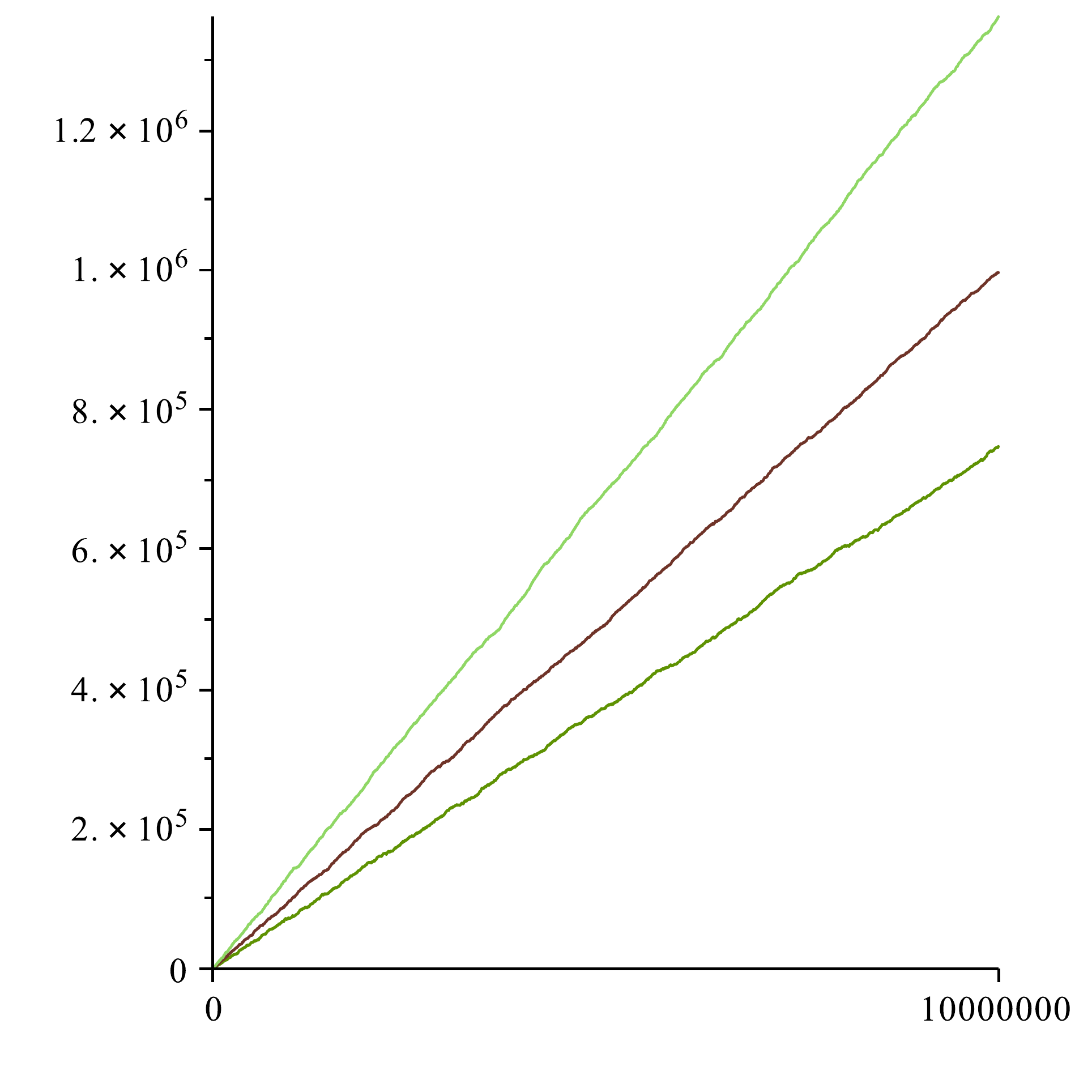}}
\put(30,0){$m=10$}
\put(145,0){\includegraphics[height=45 truemm]{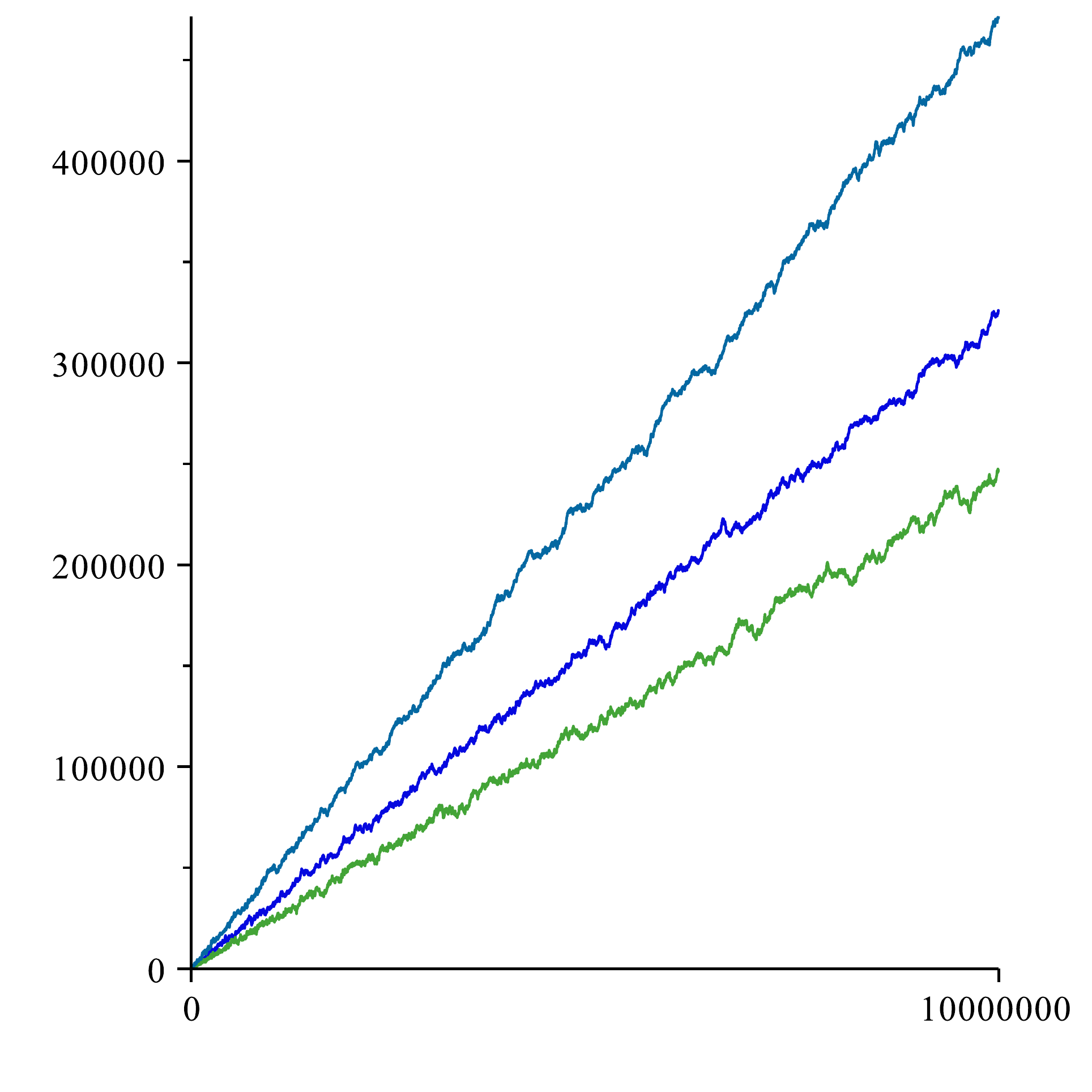}}
\put(185,0){$m=30$}
\put(300,0){\includegraphics[height=45 truemm]{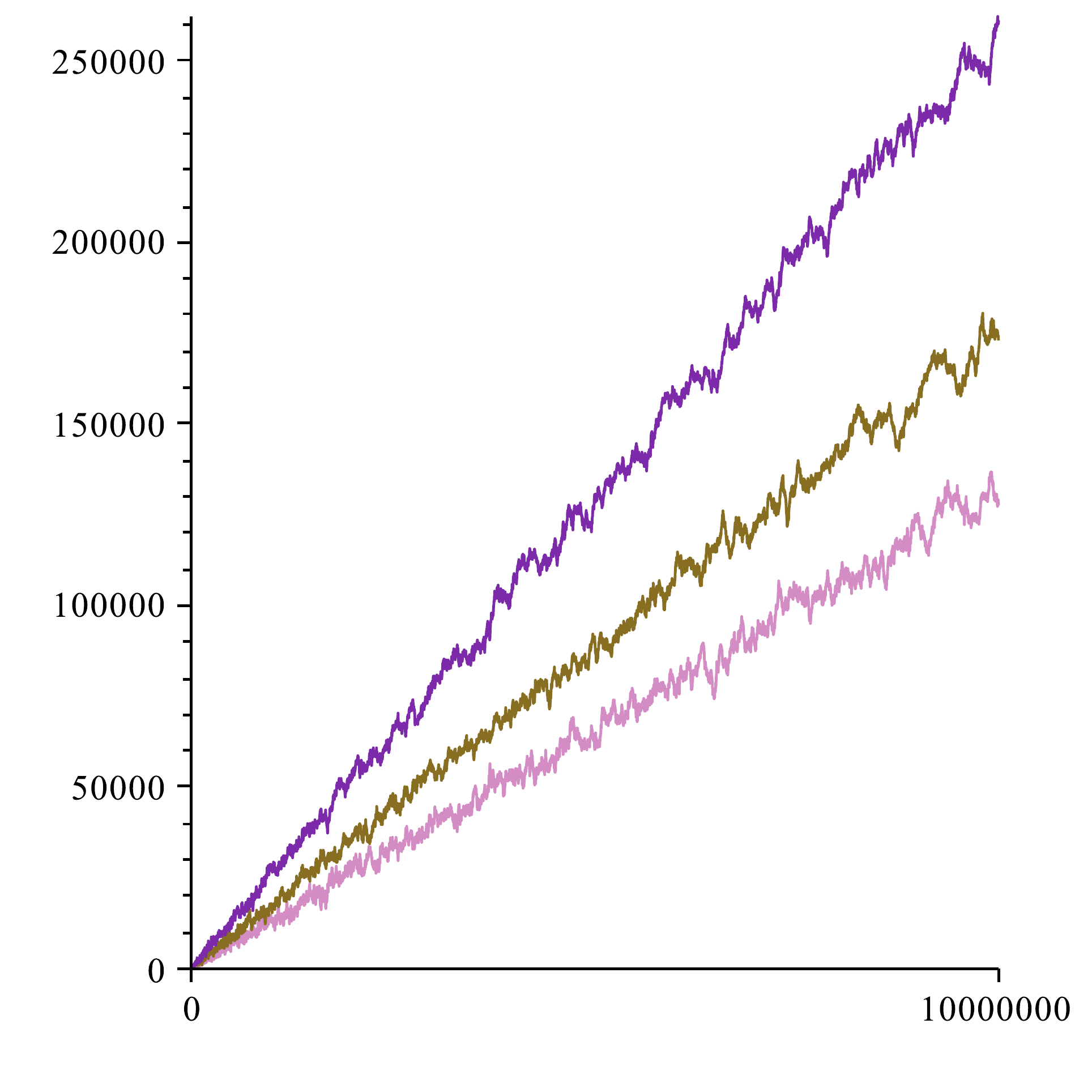}}
\put(340,0){$m=55$}
\end{picture}
\end{center}
\caption{\label{fig:mpetit} Simulations for $3$ coordinates of the gap process $\left( G_n\right) _n$ for \emph{small} $m$.}
\end{figure}

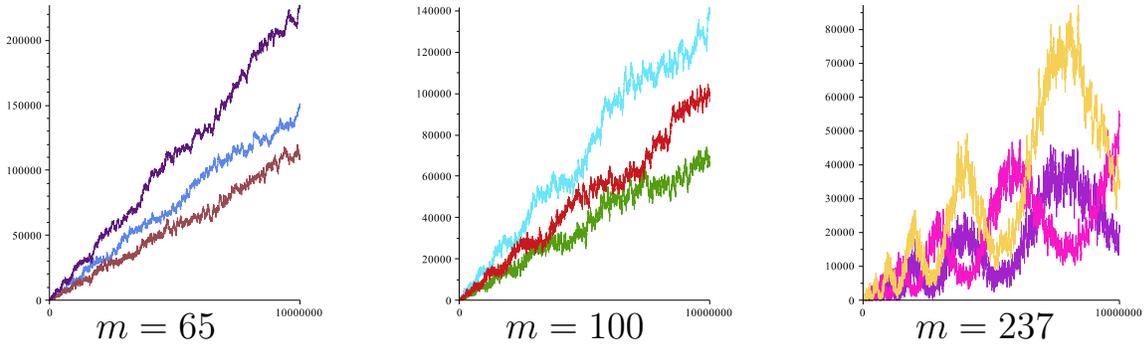
\begin{figure}[h]
\begin{center}
\begin{picture}(600,150)
\put(-10,0){\includegraphics[height=45 truemm]{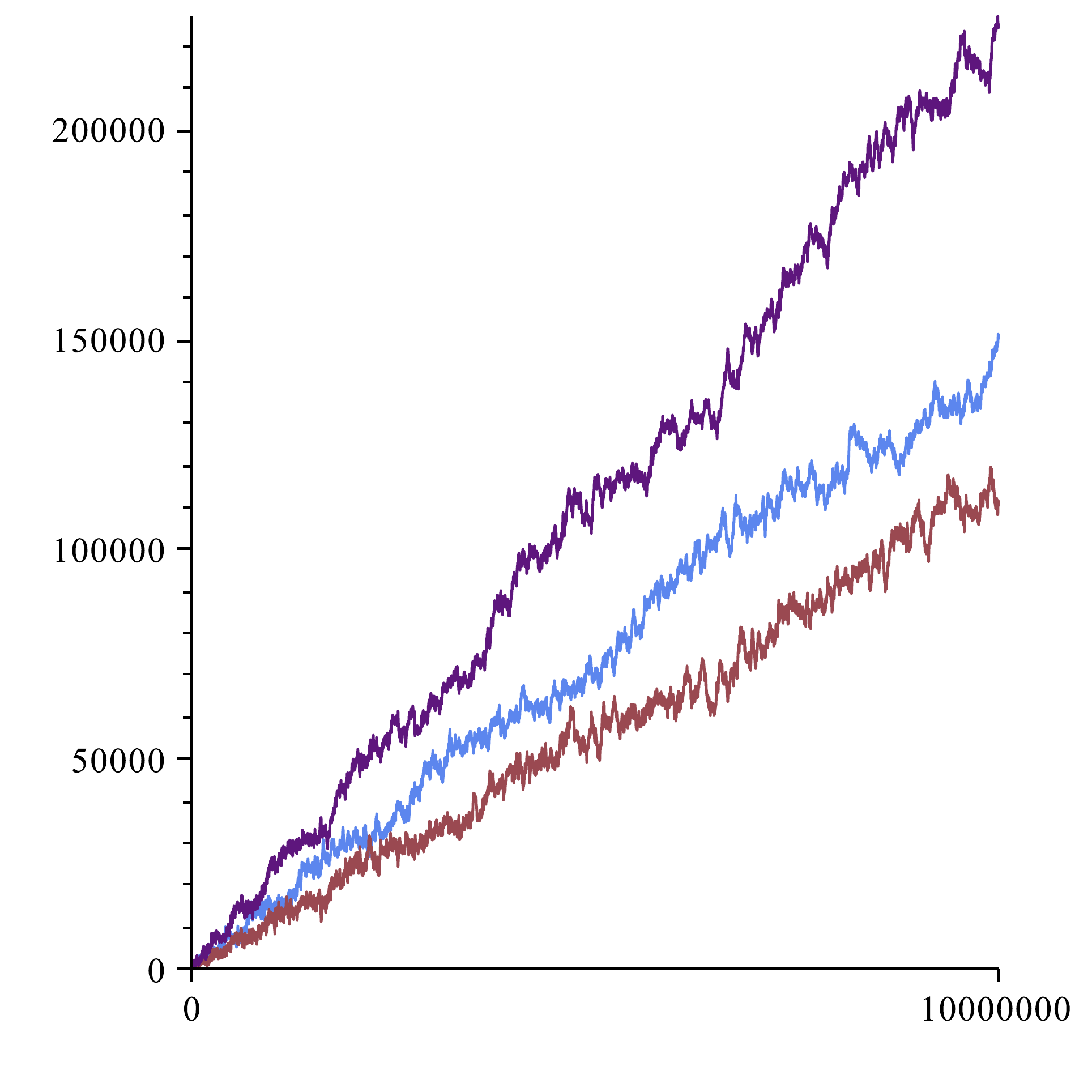}}
\put(30,0){\large$m=65$}
\put(145,0){\includegraphics[height=45 truemm]{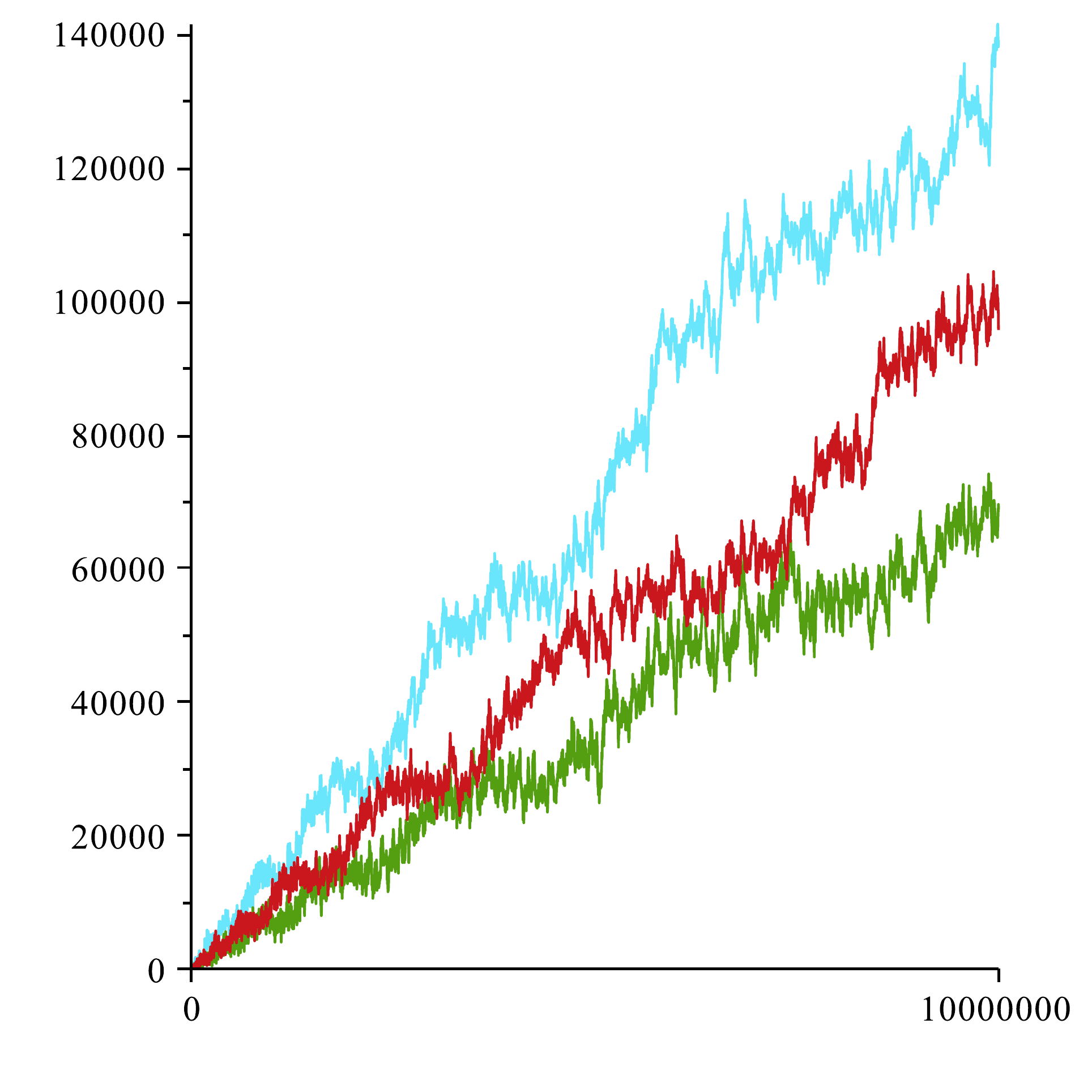}}
\put(185,0){\large$m=100$}
\put(300,0){\includegraphics[height=45 truemm]{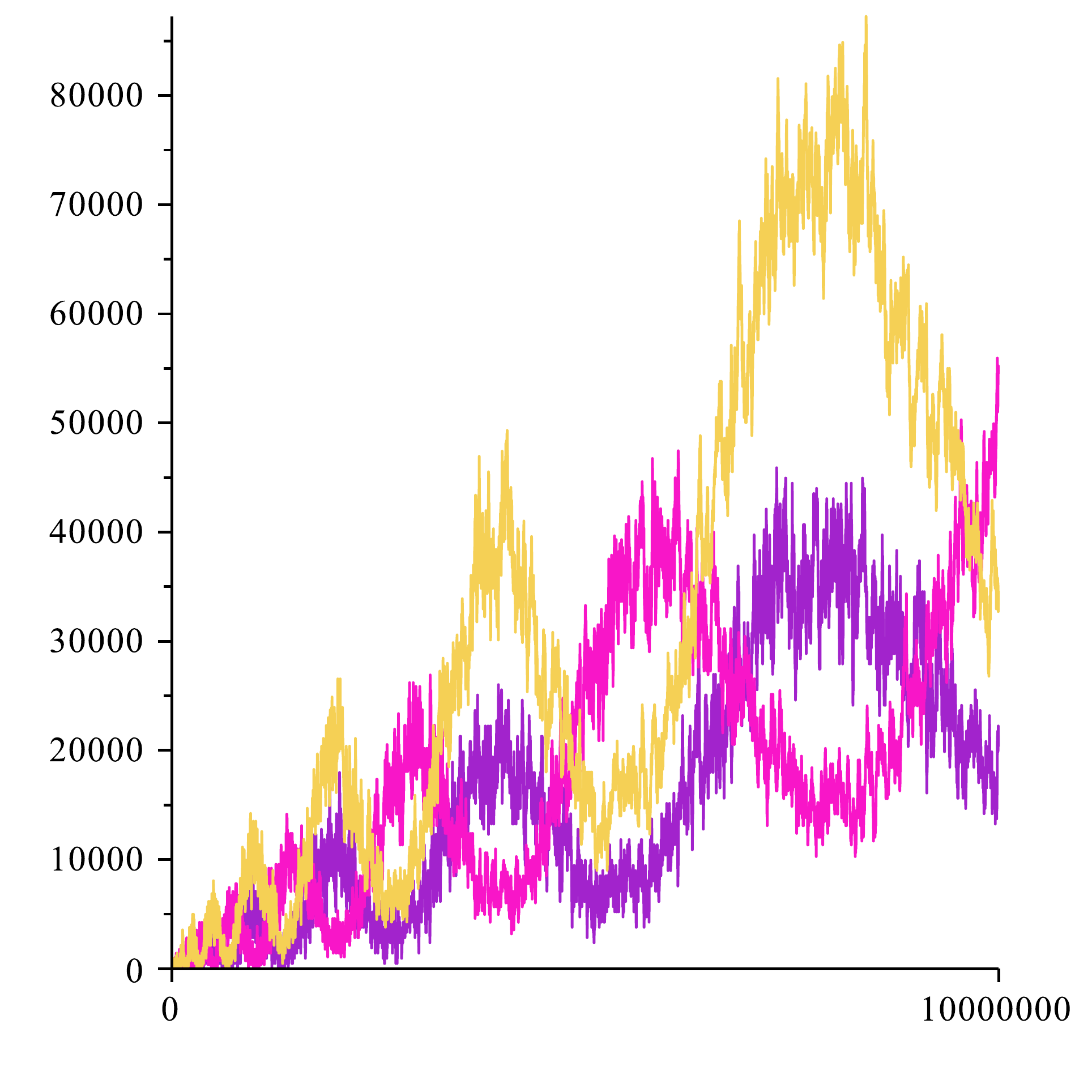}}
\put(340,0){\large$m=237$}
\end{picture}
\end{center}
\caption{\label{fig:mgrand} Simulations for $3$ coordinates of the gap process $\left( G_n\right) _n$ for \emph{large} $m$.}
\end{figure}

%\begin{itemize}
%\item[(a)]
Figures \ref{fig:mpetit} and \ref{fig:mgrand} represent the trajectories of three coordinates of the random vector~$G_n$:
for any given value $m=10$, $30$, $55$, $65$, $100$ or $237$ of the parameter, we make one random drawing of a sequence of $10^7$ keys and insert them in a B-tree.
On the pictures, the $x$-axis represents the time $n\in\{ 0,\dots ,10^7\}$ while the $y$-axis represents the number $G_n^{(k)}$ of gaps of type $k$ for $k=1$, $\lfloor m/2\rfloor$ and $m$.
%$k\in\{ 1,\lfloor m/2\rfloor ,m\}$.
In each case, the picture illustrates the almost sure asymptotics $G_n\sim nv_1$ when $n$ tends to infinity
(remember that $v_1$ is a non random $m$-dimensional vector).

\vskip 5pt
In Figure \ref{fig:mpetit}, $m$ is \emph{small} ($m=10$, $30$, $55$).
One can already catch sight of the gaussian fluctuations around the deterministic vector $nv_1$.
Notice that the variance of the gaussian limit increases with $m$, so that the amplitude of the fluctuation becomes more visible for $m=30$ and even more for $m=55$.

\vskip 5pt
In Figure \ref{fig:mgrand}, $m$ is \emph{large} ($m=65$, $100$, $237$).
On can see the almost sure oscillations around $nv_1$ appear and become more visible when $m$ grows.
Notice that they are particularly clear for $m=237$, which is the threshold value when the \emph{third} largest real part of the roots of $\chi_m$ becomes larger than $\frac 12$.
See the Appendix for more details. 

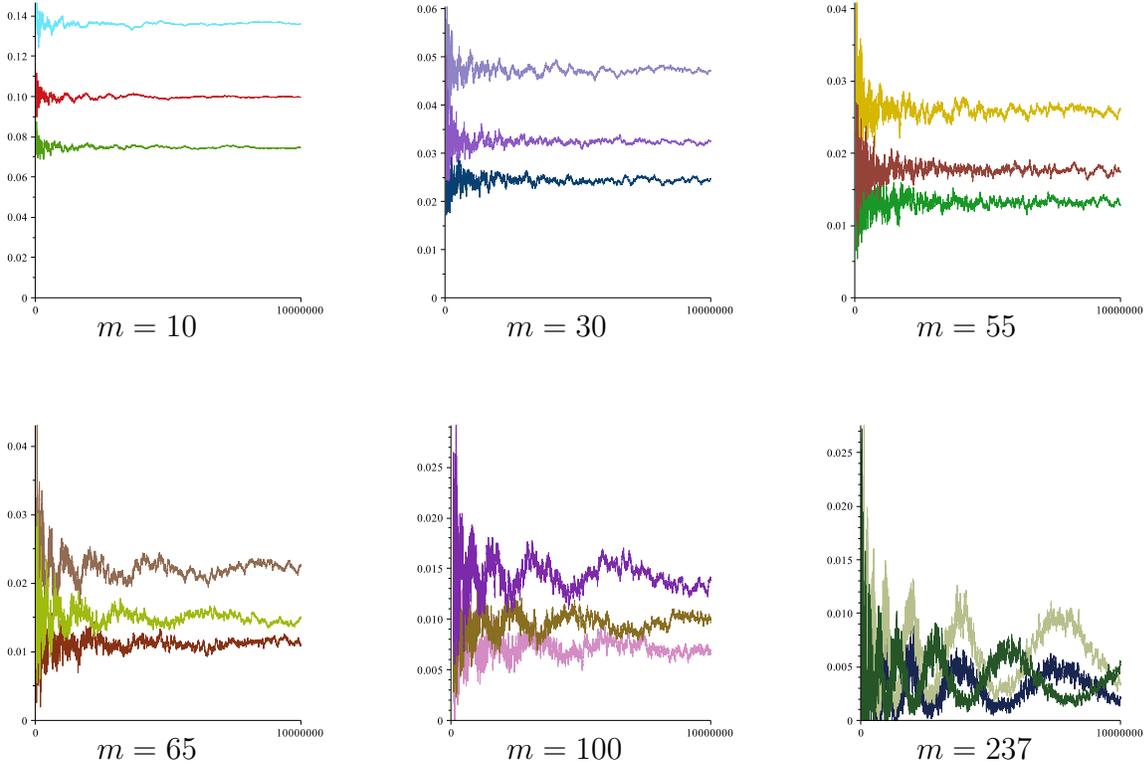
\begin{figure}[h]
\begin{center}
\begin{picture}(600,310)
\put(-10,160){\includegraphics[height=45 truemm]{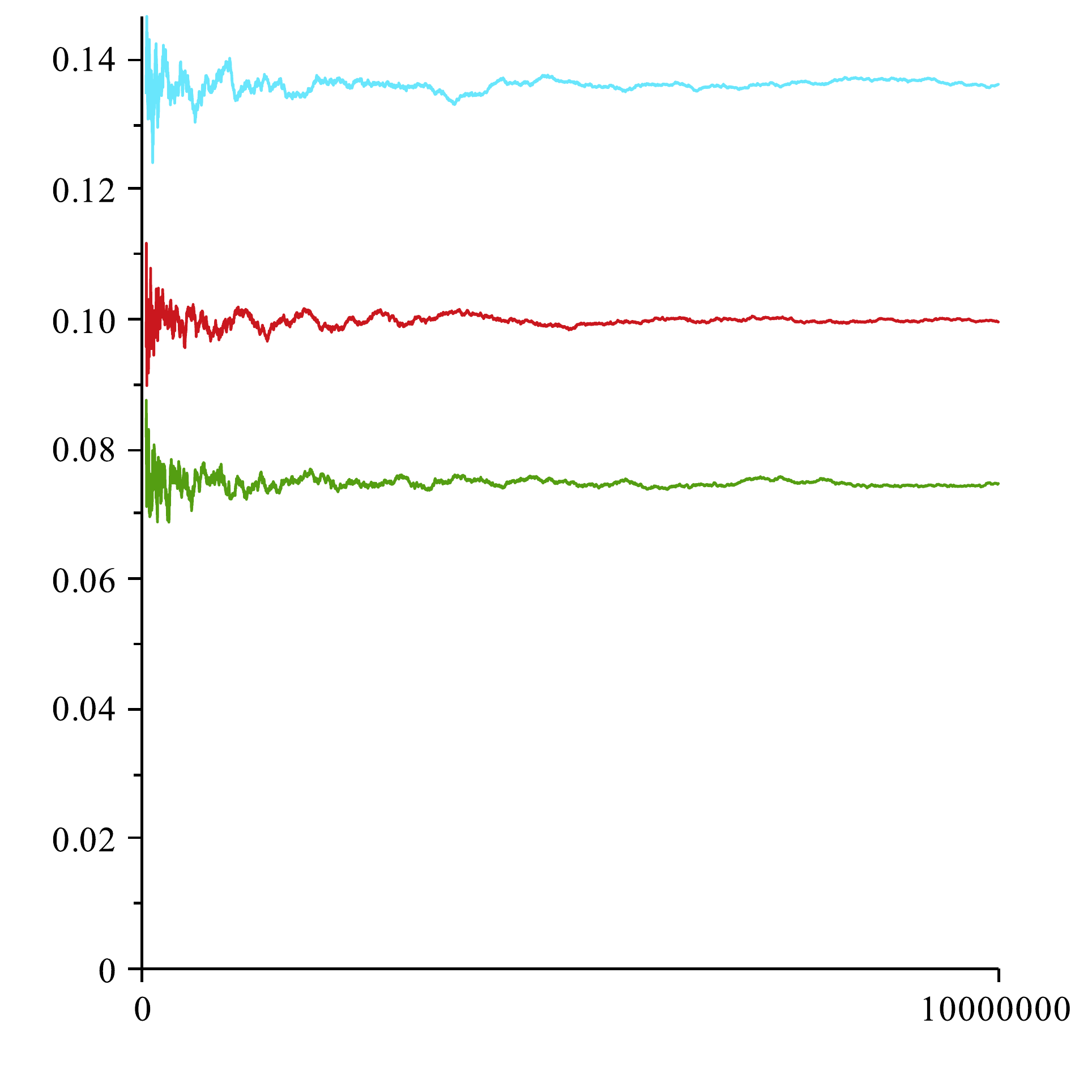}}
\put(30,160){$m=10$}
\put(145,160){\includegraphics[height=45 truemm]{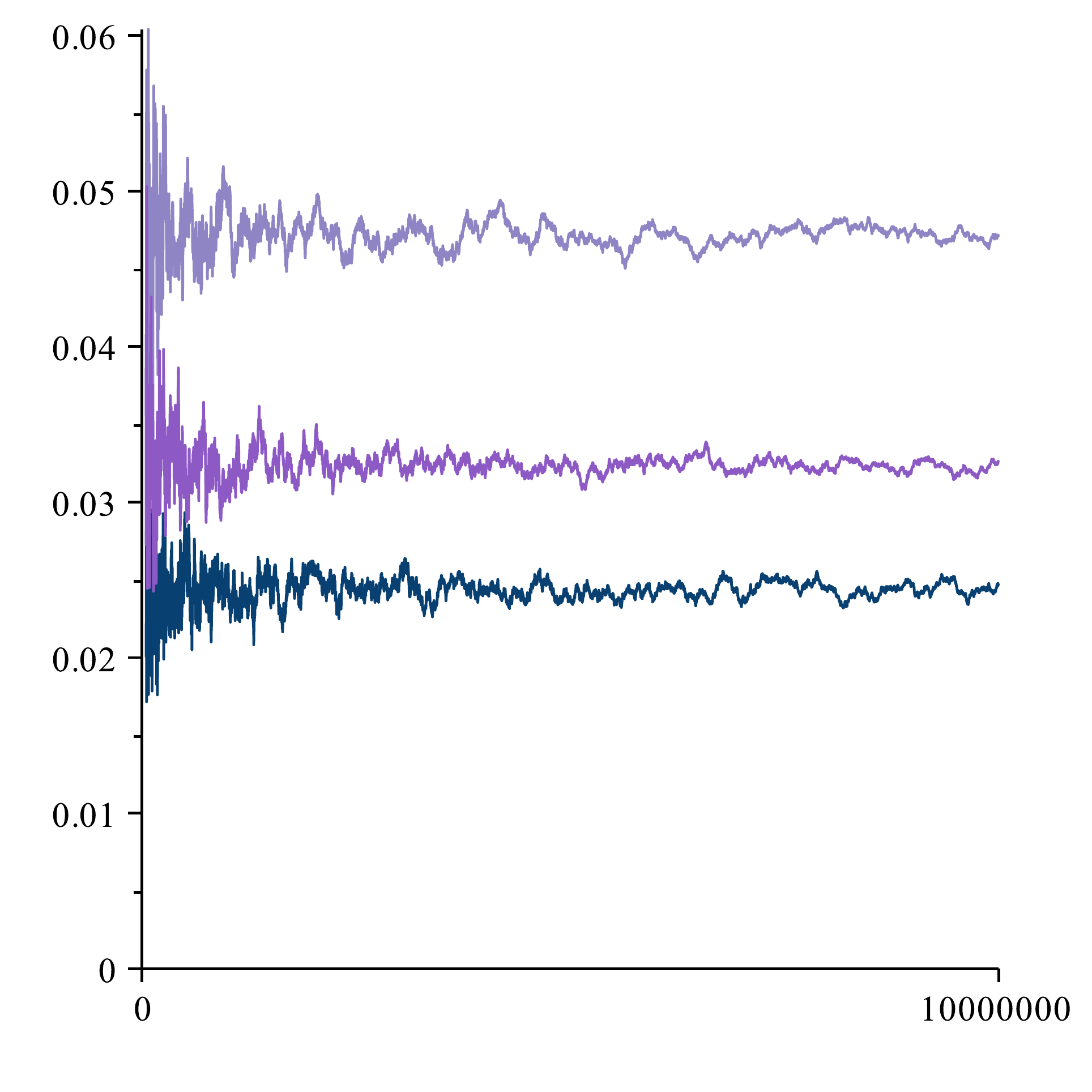}}
\put(185,160){$m=30$}
\put(300,160){\includegraphics[height=45 truemm]{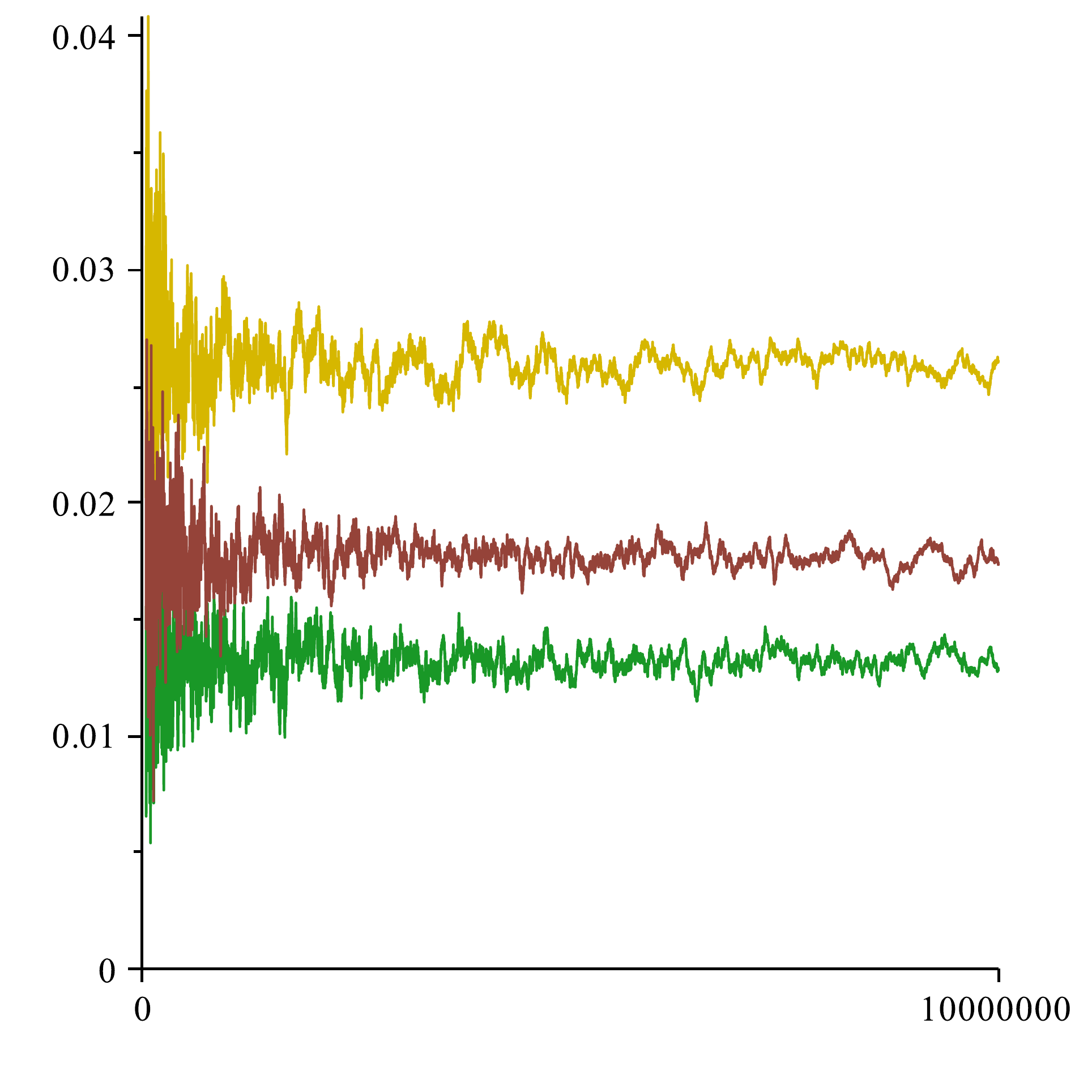}}
\put(340,160){$m=55$}
\put(-10,0){\includegraphics[height=45 truemm]{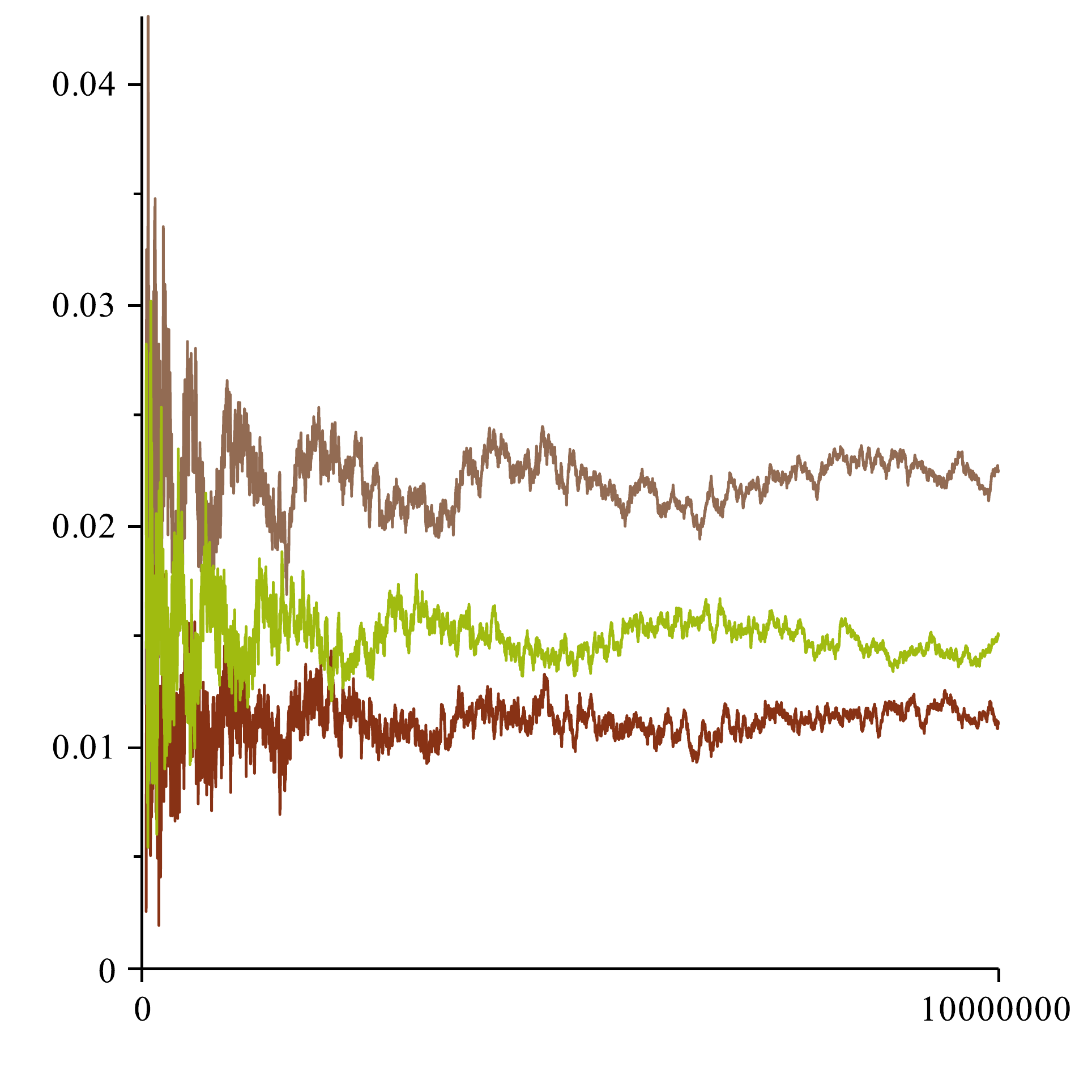}}
\put(30,0){$m=65$}
\put(145,0){\includegraphics[height=45 truemm]{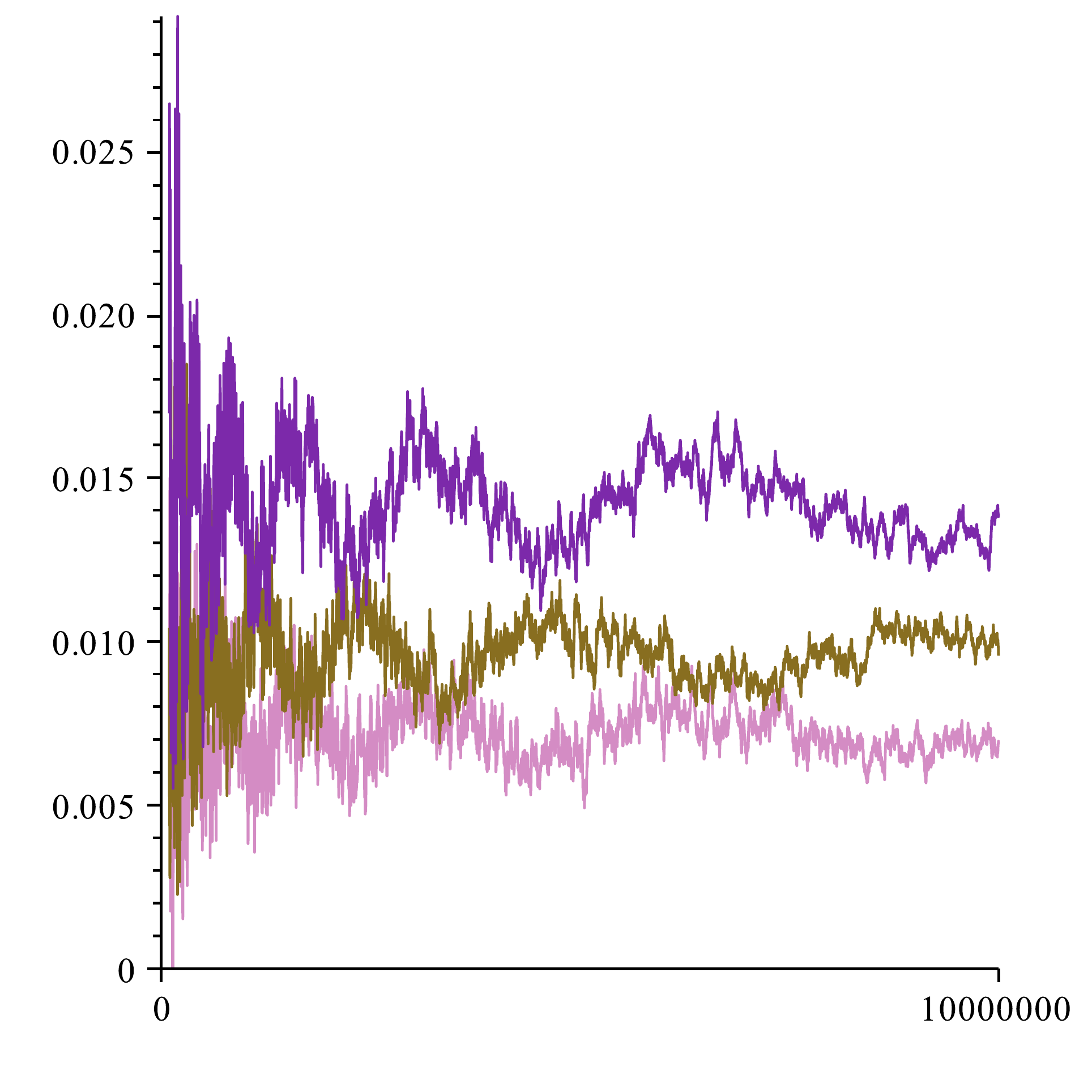}}
\put(185,0){$m=100$}
\put(300,0){\includegraphics[height=45 truemm]{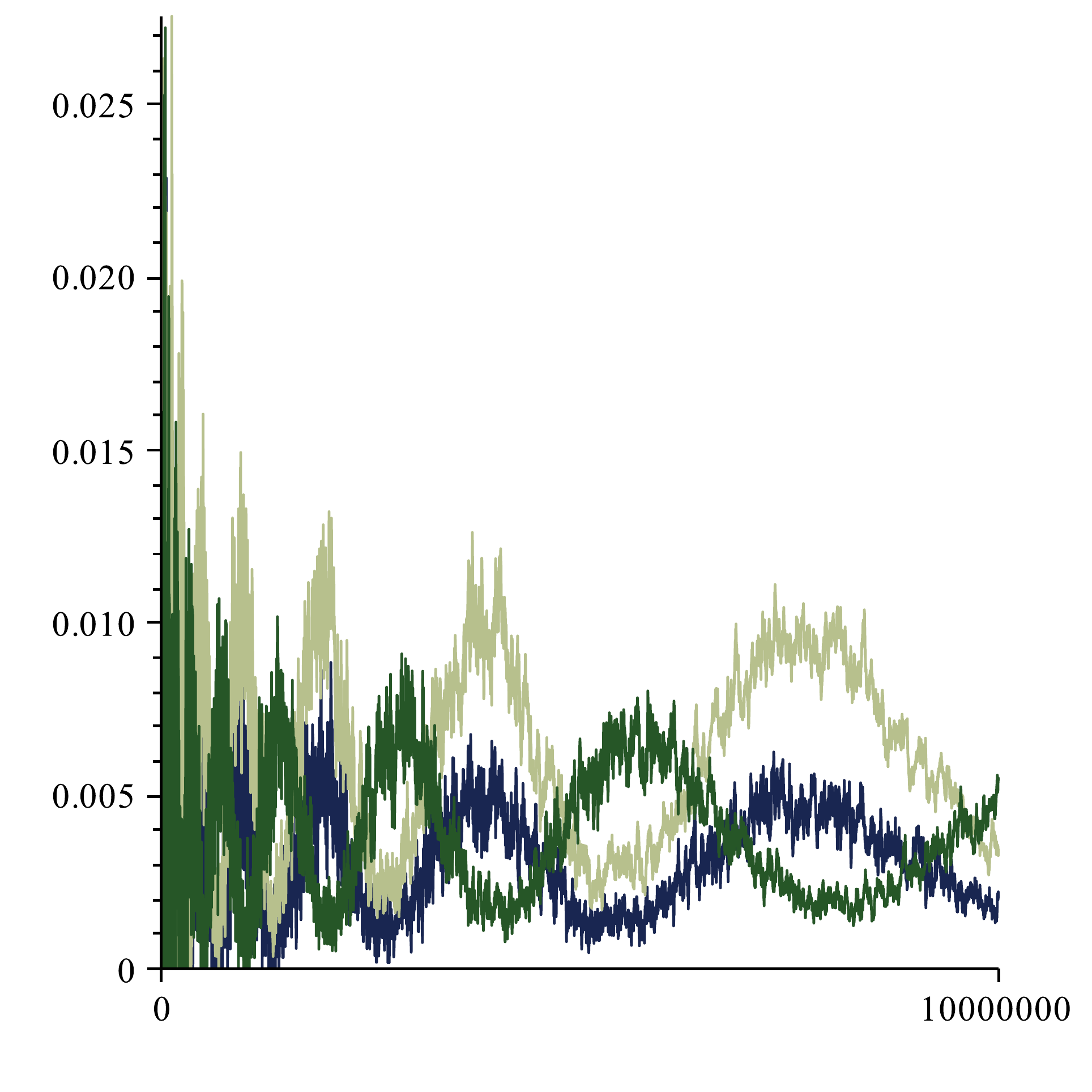}}
\put(340,0){$m=237$}
\end{picture}
\end{center}
\caption{\label{fig:deriveGn} Simulations of $3$ coordinates of  $\frac{G_n}n$ for \emph{small} and \emph{large} values of $m$.}
\end{figure}

Of course, one can make similar graphs for trajectories of the vector $\frac{G_n}n$ which converges to the deterministic vector $v_1$.
This is done in Figure~\ref{fig:deriveGn} where the convergence can be seen on the three drawn coordinates.
Once more, the fluctuations around the limit $v_1$ are of different nature depending on $m\leq 59$ or $m\geq 60$, which is also illustrated on this figure.
In particular, on can see the ``$\cos\log n$'' almost sure oscillations arise when $m\geq 60$ and become more evident when $m$ increases.

%%%%%%%%%%
\subsection{Simulations of $G_n$ after scaling}
\label{sec-simulationsGnScaling}

\begin{figure}[h]
\begin{center}
\begin{picture}(600,140)
\put(-10,0){\includegraphics[height=45 truemm]{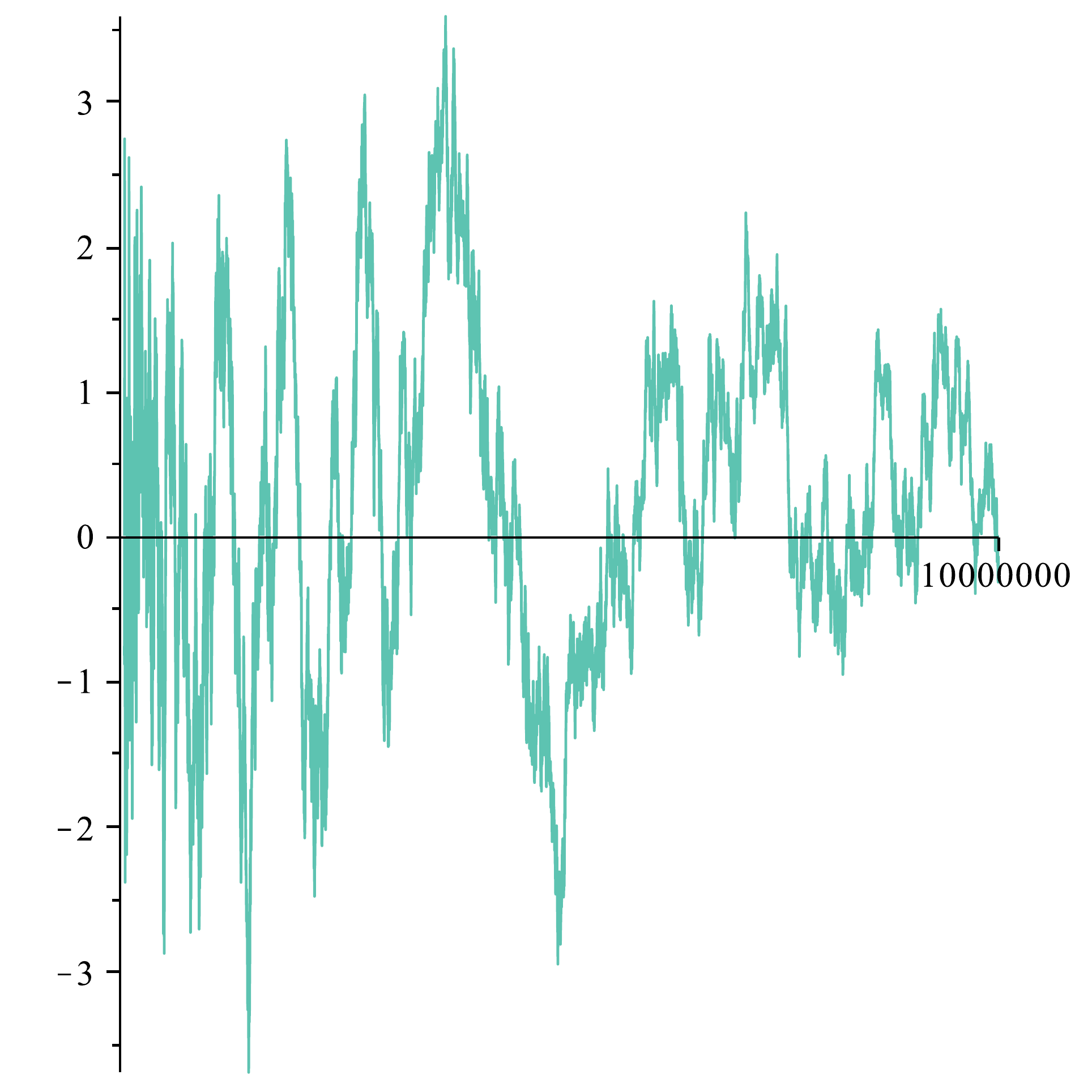}}
\put(30,-5){\large$m=10$}
\put(145,0){\includegraphics[height=45 truemm]{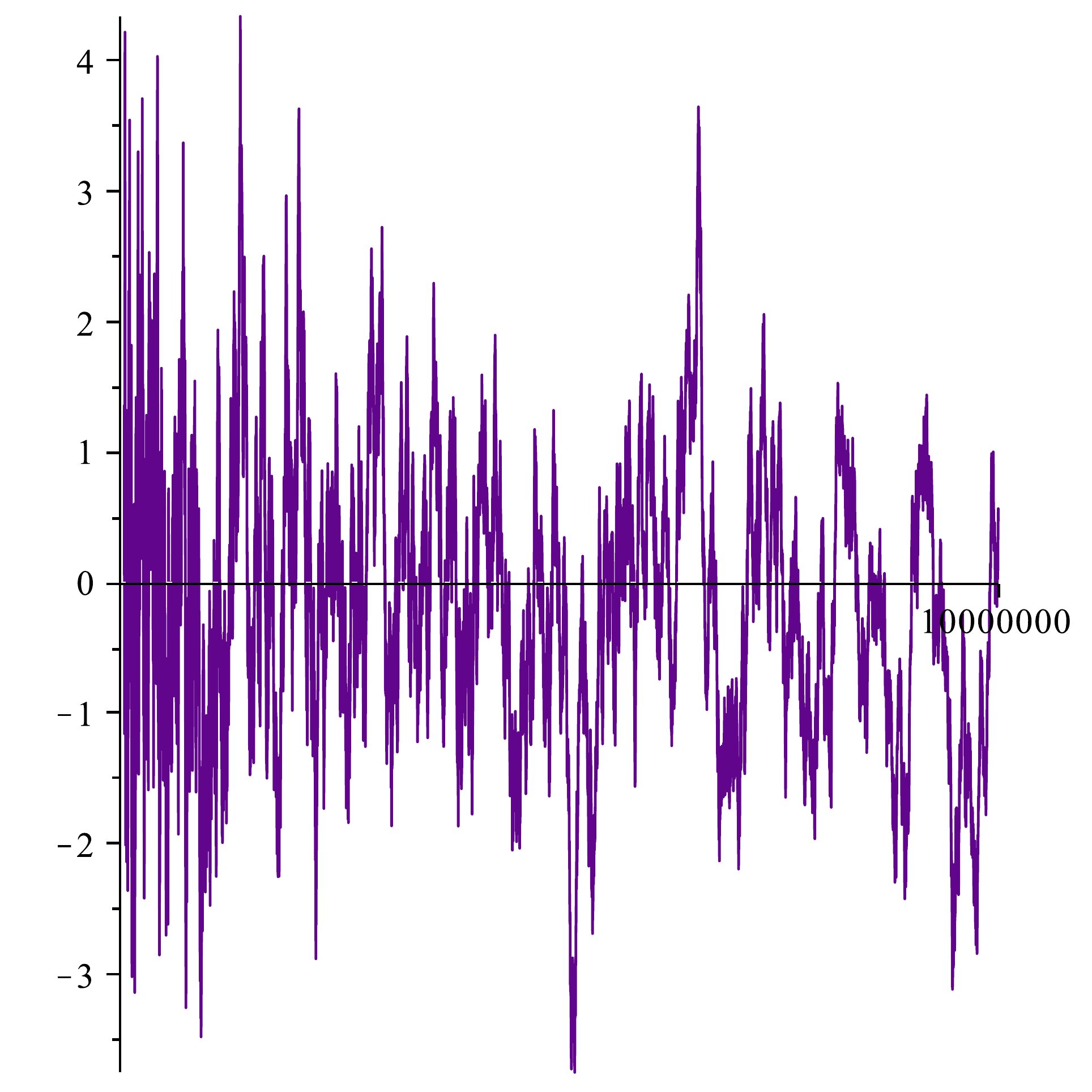}}
\put(185,-5){\large$m=30$}
\put(300,0){\includegraphics[height=45 truemm]{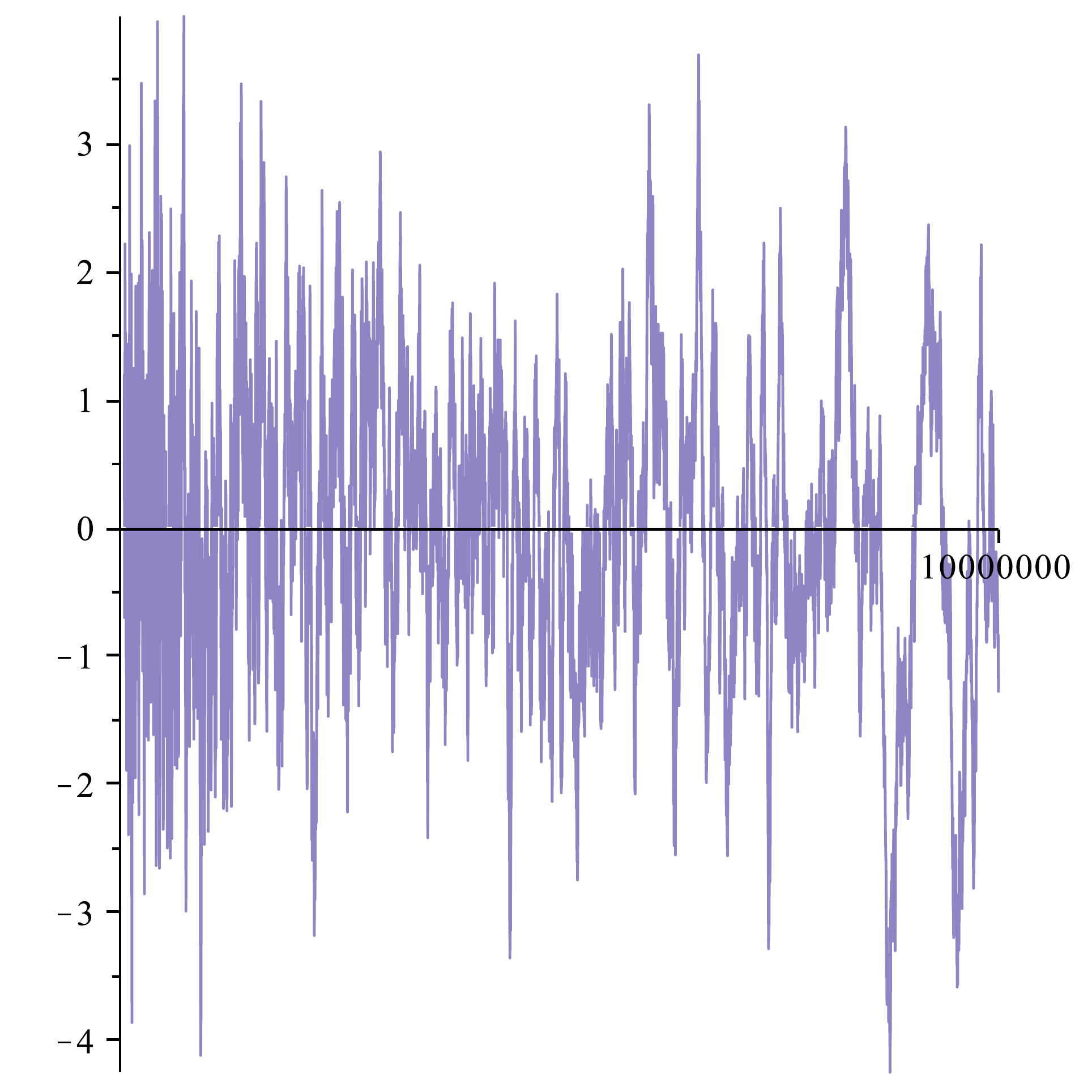}}
\put(340,-5){\large$m=55$}
\end{picture}
\end{center}
\caption{\label{fig:mpetitScaling} Simulations of one coordinate of  $G_n$ after normalisation, for \emph{small} values of $m$.}
\end{figure}
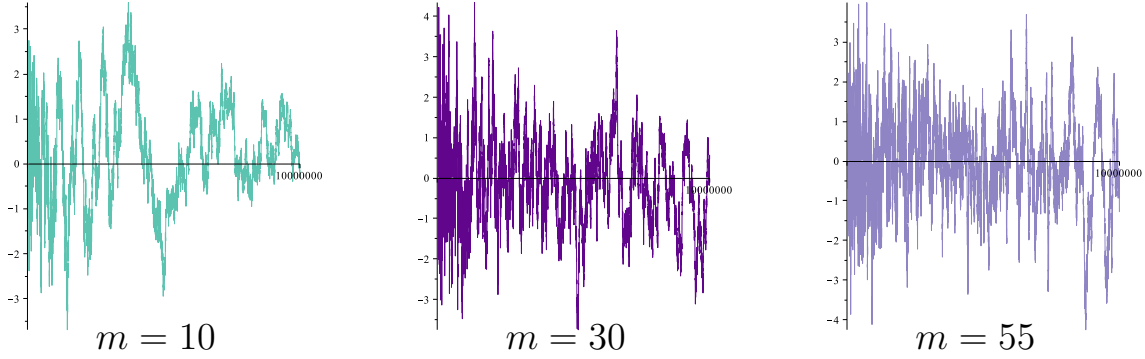

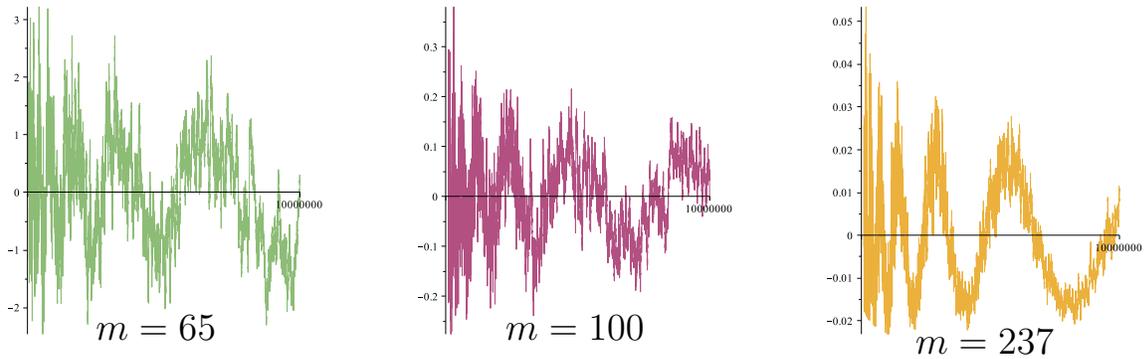
\begin{figure}[h]
\begin{center}
\begin{picture}(600,140)
\put(-10,0){\includegraphics[height=45 truemm]{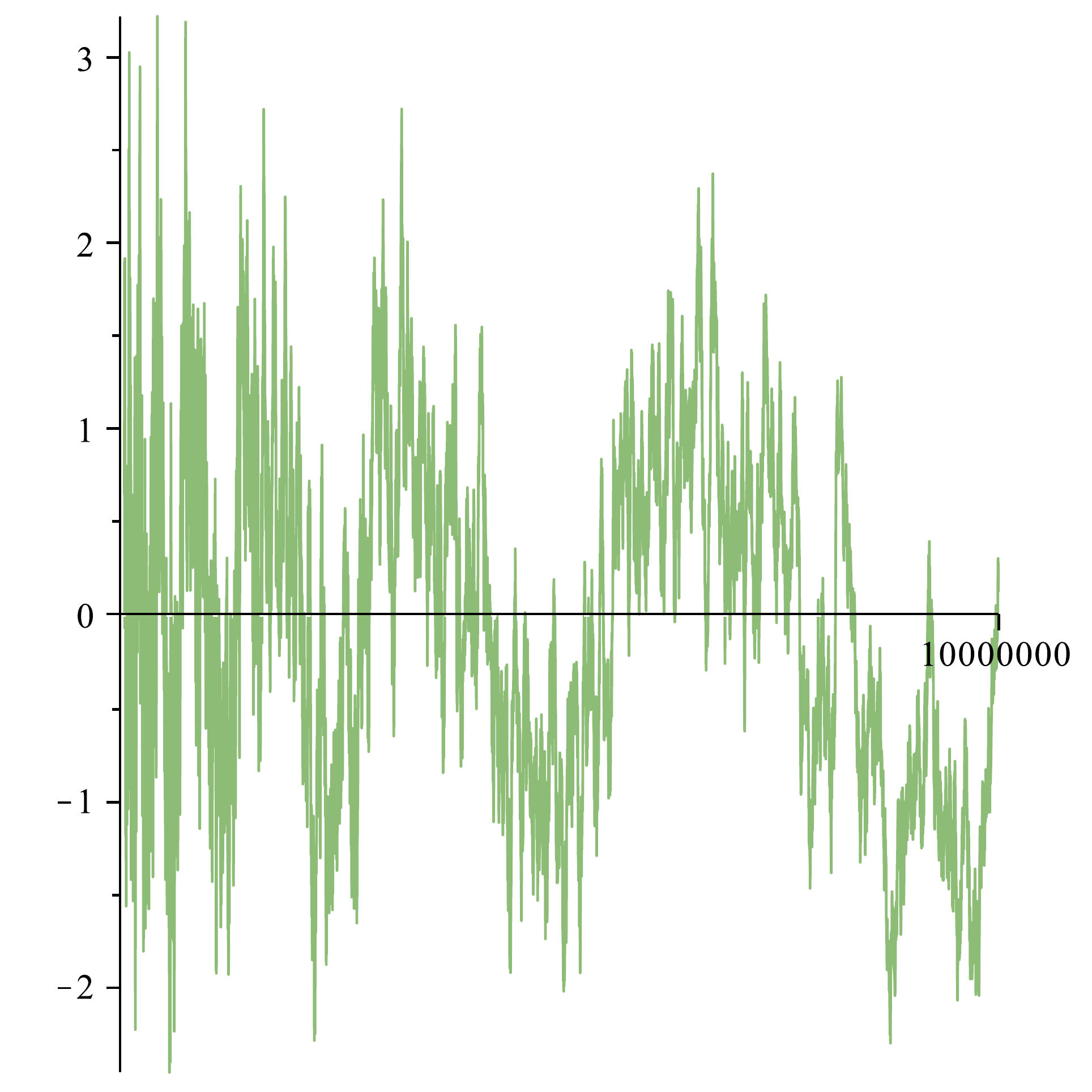}}
\put(30,0){\large$m=65$}
\put(145,0){\includegraphics[height=45 truemm]{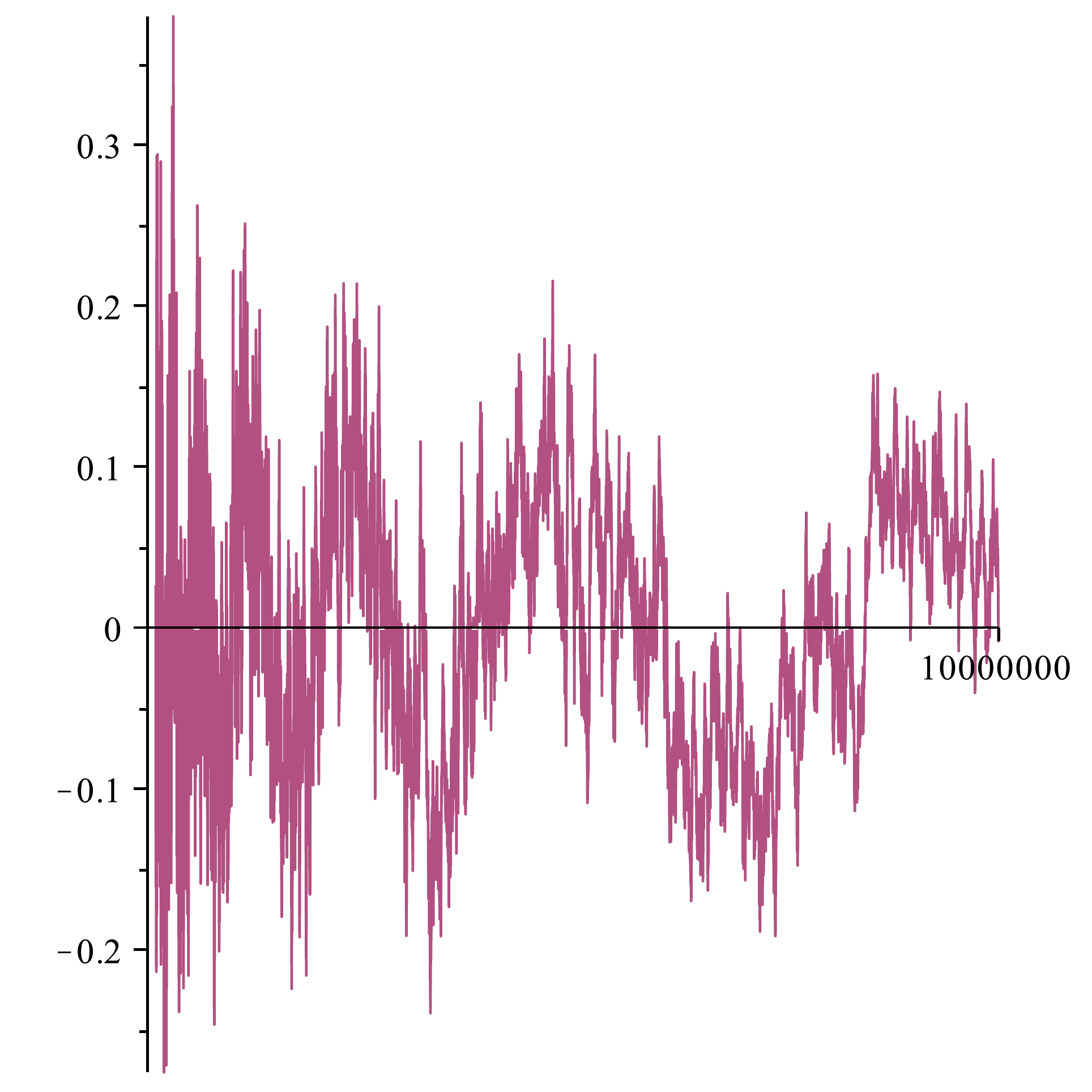}}
\put(185,0){\large$m=100$}
\put(300,0){\includegraphics[height=45 truemm]{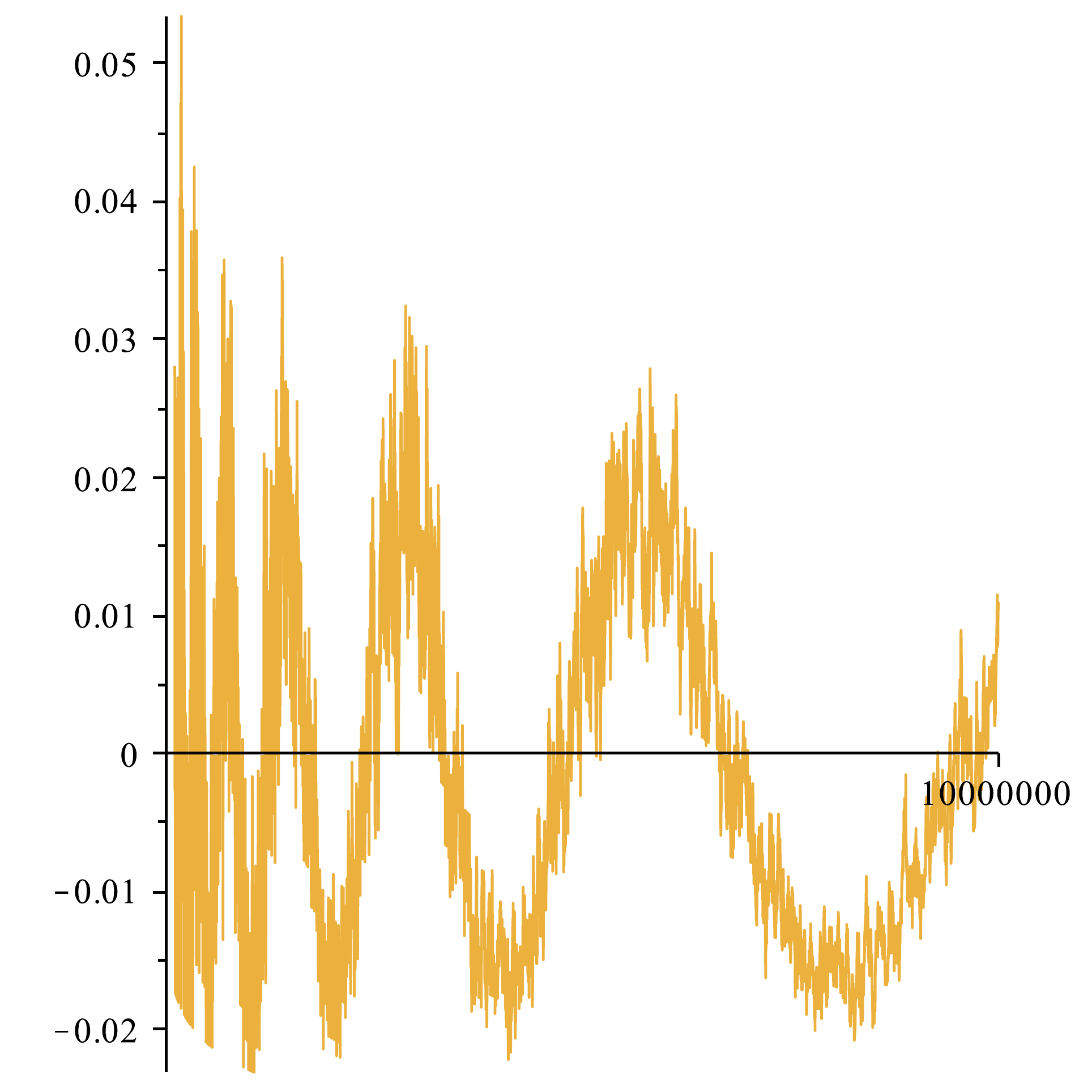}}
\put(340,-5){\large$m=237$}
\end{picture}
\end{center}
\caption{\label{fig:mgrandScaling} Simulations of one coordinate of  $G_n$ after normalization, for \emph{large} values of $m$.}
\end{figure}

A second kind of simulations focus on the possible scalings of the centered gap process $\left( G_n-nv_1\right) _n$.
In order to get convergence, according to Theorem~\ref{asymptoticsDTG}, one has to divide $G_n-nv_1$ by $\sqrt n$ when $m\leq 59$ and by $n^{\sigma _2}$ when $m\geq 60$.
Figures~\ref{fig:mpetitScaling} and~\ref{fig:mgrandScaling} represent trajectories of the median coordinate (the $\lfloor m/2\rfloor$-th) of the normalized vector process.
Hereunder, $X_n$ denotes this median coordinate $X_n=G_n^{\lfloor m/2\rfloor}$.

\vskip 10pt
Figure \ref{fig:mpetitScaling} deals with \emph{small} values of $m$, namely $m=10$, $30$, $55$ again.
On the $x$-axis, time $n\in\{ 0,\dots , 10^7\}$ ;
on the $y$-axis, the normalized coordinate $\displaystyle\frac{X_n-nv_1^{\lfloor m/2\rfloor}}{\sqrt n}$ which converges in distribution to a normal law.
Note that even if the random vector $\displaystyle\frac{G_n-nv_1}{\sqrt n}$ converges in distribution, it almost surely diverges, which is illustrated by its brownian-like trajectory.
One can refer to Gouet \cite{Gouet93} for more details on this continuous type process limit.

Figure \ref{fig:mgrandScaling} deals with $m=65$, $100$, $237$ which are \emph{large} values of $m$.
On the $y$-axis: the normalized coordinate $\displaystyle\frac{X_n-nv_1^{\lfloor m/2\rfloor}}{n^{\sigma _2}}$, which is almost surely equivalent to some $\rho\cos\left( \tau _2\log n+\varphi\right)$ when $n$ tends to infinity, where $\rho$ is a positive random variable (random amplitude), $\varphi$ a $[0,2\pi [$-valued random variable (random phase) and $\tau _2$ the imaginary part of the complex eigenvalue $\lambda _2=\sigma _2+i\tau _2$.
The random variables $\rho$ and $\varphi$ are proportional to the module and the argument of the complex-valued random variable $W_m$ in Theorem~\ref{asymptoticsDTG}.

%%%%%%%%%%%%%%%%%%%%%%%%%%%%%%%%%%
\section{Appendix. The phase transition and $\sigma _2(m)$}
\label{appendix}

The phase transition that occurs for B-trees with parameter $m$ relies on the roots of the characteristic polynomial
$$
\chi _m(X)=
\prod _{k=m}^{2m-1}(X+k)-\frac{(2m)!}{m!}.
%=\frac{(2m)!}{m!}\left( \prod _{k=m}^{2m-1}\frac{X+k}{1+k}-1\right)
$$
Denote by $\lambda _2=\lambda _2(m)$ the root of $\chi _m$ having the second largest real part and a positive imaginary part.
Denote also $\sigma _2=\sigma _2(m)$ the real part of $\lambda _2(m)$.

As shown in Section~\ref{sec-transition}, the B-tree admits a Gaussian central limit theorem when $m\leq 59$ (small r\'egime) whereas it admits an almost sure nonnormal fluctuation term of order $n^{\sigma _2(m)}$ around the drift when $m\geq 60$ (large r\'egime).
Coming from P\'olya urn theory, this asymptotic behaviour depends on whether $\sigma _2(m)<1/2$ (small r\'egime) or $\sigma _2(m)>1/2$ (large r\'egime).

\vskip 10pt
%\subsection{Expansion of $\sigma _2$ and phase transition}

\begin{center}
\begin{figure}[h!]
\begin{center}
\includegraphics[height=210pt]{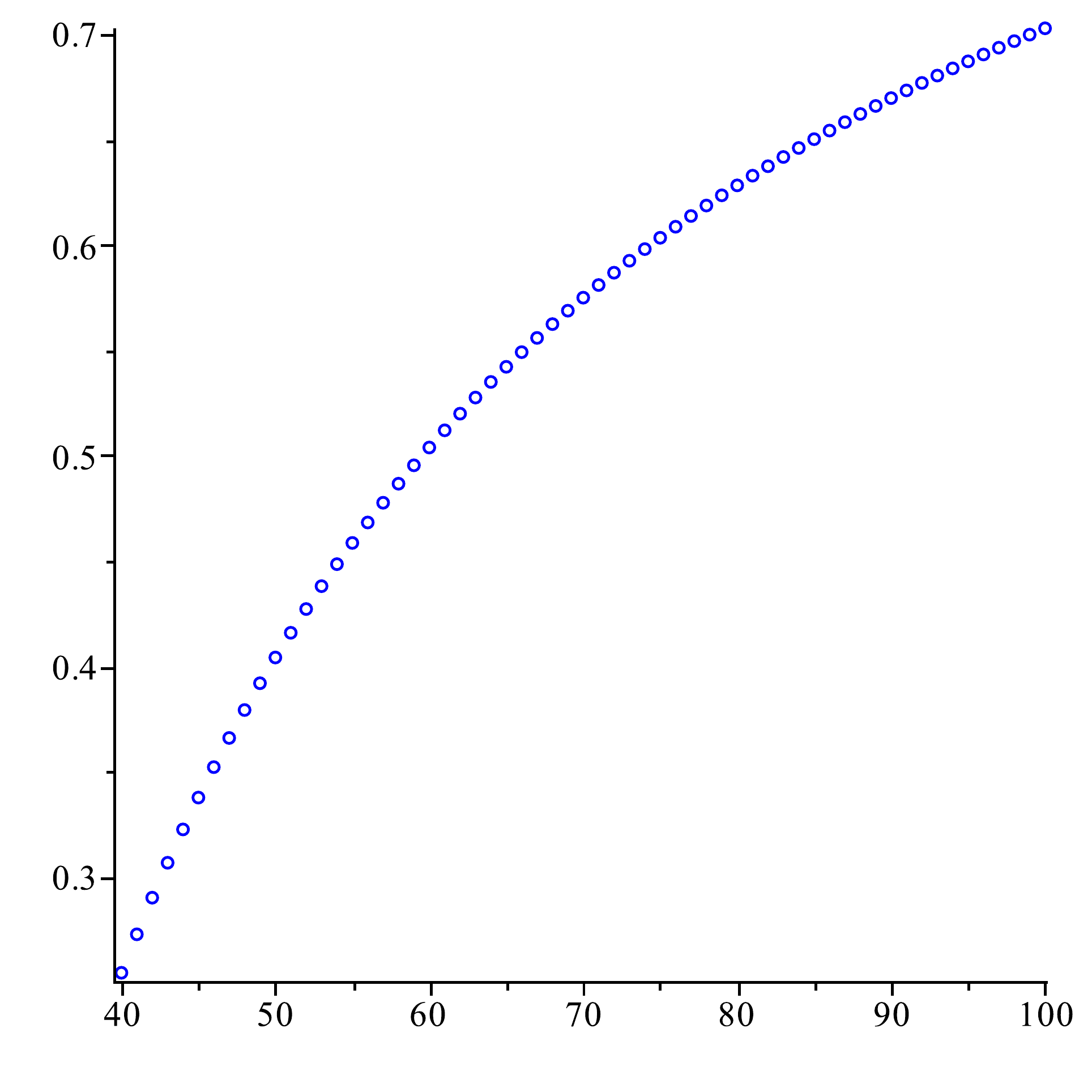}
\includegraphics[height=210pt]{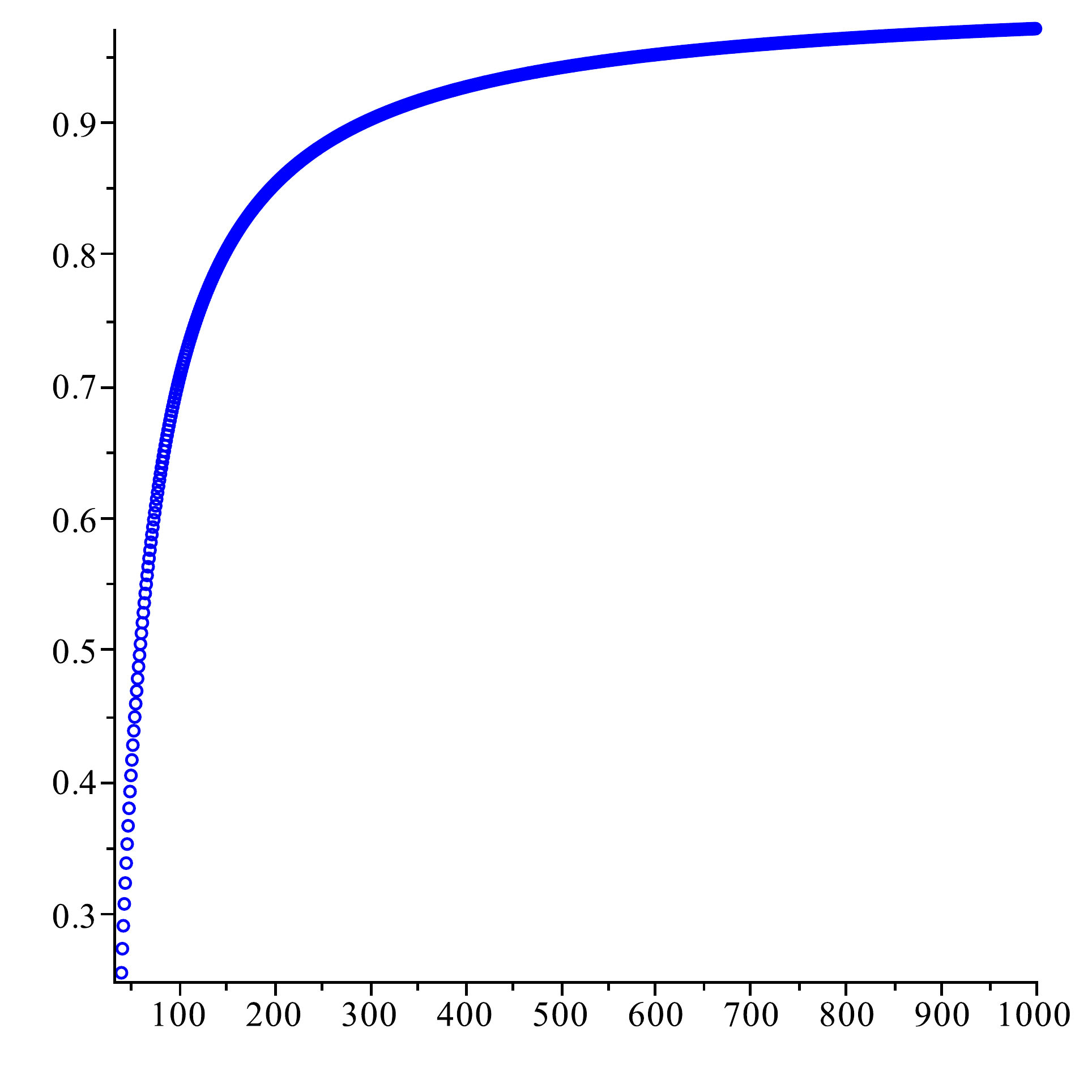}
\end{center}
\caption{%
\label{fig-sigma}
The graph of the sequence $m\mapsto\sigma _2(m)$.}
\end{figure}
\end{center}

Let $F$ the two-variable meromorphic function defined by
$$
F(x,y)=
\frac
{\Gamma\left(x+2y\right)\Gamma\left( 1+y\right)}
{\Gamma\left( 1+2y\right)\Gamma\left( x+y\right)}
$$
where $\Gamma$ denotes Euler's Gamma function.
For a given $m\geq 2$, $\lambda _2=\lambda _2(m)=\sigma _2+i\tau _2$ is the root of equation $F\left( X,m\right) =1$ having the the second largest real part $\sigma _2$ (the first one being reached by the evident root~$1$) and a positive imaginary part~$\tau _2$.
Denote by $\psi$ the classical Digamma function, the logarithmic derivative of Euler's Gamma.
Since $\frac{\partial}{\partial x}F(x,1/y)=\psi\left( x+2/y\right) -\psi\left(x+1/y\right) =\log 2+O(y)$ as $y$ tends to $0$, the analytic implicit function theorem shows that $\lambda _2(m)$ is an analytic function of $1/m$ as $m$ tends to $+\infty$.
Using the expansion
$$
\log\Gamma (z)=z\log z-z-\frac 12\log z+\frac 12\log 2\pi+O\left(\frac 1z\right) ~{\rm mod}~2i\pi
$$
as $|z|$ tends to infinity (the \emph{mod} coming from the determination of the logarithm), writing  $\lambda _2(m)$ as a power series in $1/m$ and  putting the first terms of this expansion of in the equation $\log F\left(\lambda _2(m),m\right)=0~{\rm mod~2i\pi}$, one gets by identification
$$
\left\{
\begin{array}{l}
\displaystyle
\sigma _2(m)=1-\frac{\pi ^2}{\log ^32}\times\frac 1m+O\left(\frac 1{m^2}\right)
\\ \\
\displaystyle
\tau _2(m)=\frac{2\pi}{\log 2}+\frac{\pi}{2\log ^22}\times\frac 1m+O\left(\frac 1{m^2}\right)
\end{array}
\right.
$$
as $m$ tends to infinity.
The graph of the function $m\mapsto\sigma _2(m)$ is given in Figure~Ê\ref{fig-sigma}.
Numerical values of the expansions give $\sigma _2\approx 1-29.63/m+\dots$ while $\tau _2\approx 9.06+3.27/m+\dots$

Moreover, the numerical values of $\sigma _2$ around $m=60$ are the following ones, showing more accurately that $\sigma _2(m)<1/2$ if, and only if $m\leq 59$.
These numerical values have been computed applying the Newton method to the $\chi _m$ function starting from the point $0.5+9.0i$, as suggested by the above expansions of $\sigma _2$ and $\tau _2$.

\begin{center}
\begin{tabular}{|c|c|}
\hline
$m$&$\sigma _2(m)$\\
\hline
57&0.4775726941\\
\hline
58&0.4866133472\\
\hline
59&0.4953467200\\
\hline
60&0.5037882018\\
\hline
61&0.5119521623\\
\hline
62&0.5198520971\\
\hline
\end{tabular}
\end{center}

In order to justify the choice of $m=237$ in our drawings, denote by $\sigma _3(m)$ the third largest real part of the roots of $\chi _m$.
 The threshold value when $\sigma _3(m)$ becomes larger than $\frac 12$ is $m=237$.
 Using general statements on P\'olya urns, this shows that a second almost sure phenomenon with magnitude $n^{\lambda _3}$ is added to the one we describe in Theorem~\ref{asymptoticsDTG} as soon as $m\geq 238$.
 That is the reason why the above figures have been selected for $m=237$; indeed, for $m$ increasing from $60$ until $237$, the asymptotic expansion of $G_n$ contains the oscillating term of amplitude $n^{\sigma _2(m)}$ more and more visible compared to brownian terms in $n^{\frac 12}$. For $m>237$, the second oscillating term of amplitude $n^{\sigma _3(m)}$ appears, making the $n^{\sigma _2(m)}$ oscillation less visible.
The numerical values of $\sigma _3(m)$ around $m=237$ are the following ones.

\vskip 20pt

\begin{center}
\begin{tabular}{|c|c|}
\hline
$m$&$\sigma _3(m)$\\
\hline
236&0.4971039325\\
\hline
237&0. 4992277960\\
\hline
238&0.5013338161\\
\hline
239&0.5034221856\\
\hline
\end{tabular}
\end{center}

\vskip 30pt
{\large{\bf Acknowledgements}
}
The authors kindly thank Karine Zeitouni for valuable discussions on B-trees and Andrea Sportiello for insightful comments.
%%%%%%%%%%%%%%%%%%%%%%%%%%%%%%%%%%%%%%%%%%%%%
\bibliographystyle{plain}
\bibliography{biblioChaGarPouTon}
%%%%%%%%%%%%%%%%%%%%%%%%%%%%%%%%%%%%%%%%%%%%

\end{document}